\documentclass[a4paper,11pt]{article}

\usepackage{amsmath}
\usepackage{amssymb}
\usepackage{pgfplots}
\pgfplotsset{compat=newest}
\usepackage{amsthm}
\usepackage{mathtools}
\usepackage{booktabs}
\usepackage{enumerate}
\usepackage{cleveref}
\usepackage{geometry}
\usepackage{graphicx}
\usepackage{subcaption}
\usepackage{cancel}
\usetikzlibrary{shapes,snakes}
\usepackage{3dplot} 

\newgeometry{vmargin={35mm}, hmargin={20mm,20mm}}   

\newtheorem{theorem}{Theorem}

\newtheorem{problem}[theorem]{Problem}

\theoremstyle{definition}
\newtheorem{definition}[theorem]{Definition}
\newtheorem{assumption}[theorem]{Assumption}

\DeclareMathOperator{\rank}{rank}

\DeclareMathOperator{\Ima}{Im}

\newcommand{\R}{\mathbb{R}}
\newcommand{\SH}[2]{Y_{#1}^{#2}}

\newcommand{\SHn}[1]{\Upsilon_{\text{#1}}}
\newcommand{\SHnB}[1]{\bar{\Upsilon}_{\text{#1}}}
\newcommand{\SHne}[2]{\Upsilon_{\text{#1}}^\text{e,#2}}
\newcommand{\SHno}[2]{\Upsilon_{\text{#1}}^\text{o,#2}}
\newcommand{\SHN}[1]{\Upsilon_{\text{#1}}}
\newcommand{\SHNe}[2]{\Upsilon_{\text{#1}}^\text{e,#2}}
\newcommand{\SHNo}[2]{\Upsilon_{\text{#1}}^\text{o,#2}}
\let\oldint\int
\newcommand{\eqdef}{\mathrel{\mathop:}=}
\renewcommand\int{\oldint\limits}
\newcommand{\N}{\mathbb{N}}

\newcommand{\spanVecSpace}[1]{\text{span}\lbrace #1 \rbrace}
\newcommand{\domain}{G}
\newcommand{\PN}[1]{\text{$P_{\text{#1}}$}}
\newcommand{\normpdf}[1]{e^{- \left(\frac{#1-\mu_{#1}}{\sqrt{2} \sigma_{#1}} \right)^2}}
\newcommand{\US}[1]{\mathbb{S}^{#1}}

\newcommand{\unp}{u^{n_{+}}}
\newcommand{\unm}{u^{n_{-}}}
\newcommand{\Qnp}{Q^{n_{+}}}
\newcommand{\Qnm}{Q^{n_{-}}}

\newenvironment{nalign}{
	\begin{equation}
	\begin{aligned}
		}{
	\end{aligned}
	\end{equation}
	\ignorespacesafterend
}

\date{}
\title{Stable Boundary Conditions and Discretization for\\ \PN{N} Equations}
\author{Jonas B\"unger, Neeraj Sarna and Manuel Torrilhon \\[1ex]
	{\small Center for Computational Engineering \& Department of Mathematics}\\
	\small RWTH Aachen University, Germany\footnote{Mathematics (CCES), Schinkelstr. 2, 52062 Aachen, Germany}\\
	(\texttt{\small buenger@mathcces.rwth-aachen.de,sarna@mathcces.rwth-aachen.de,mt@mathcces.rwth-aachen.de})}

\begin{document}
	\maketitle
	\begin{abstract}
	A solution to the linear Boltzmann equation satisfies an energy bound, which reflects a natural fact: The energy of particles in a finite volume is bounded in time by the energy of particles initially occupying the volume augmented by the energy transported into the volume by particles entering the volume over time. In this paper, we present boundary conditions (BCs) for the spherical harmonic (\PN{N}) approximation, which ensure that this fundamental energy bound is satisfied by the \PN{N} approximation. Our BCs are compatible with the characteristic waves of \PN{N} equations and determine the incoming waves uniquely. Both, energy bound and compatibility, are shown on abstract formulations of \PN{N} equations and BCs to isolate the necessary structures and properties.
	The BCs are derived from a Marshak type formulation of BC and base on a non-classical even/odd-classification of spherical harmonic functions and a stabilization step, which is similar to the truncation of the series expansion in the \PN{N} method.
	We show that summation by parts (SBP) finite differences on staggered grids in space and the method of simultaneous approximation terms (SAT) allows to maintain the energy bound also on the semi-discrete level.
\end{abstract}

	\section{Linear Boltzmann and Radiative Transfer Equation}
The transport of particles in a background medium is governed by the \emph{linear Boltzmann equation}\cite{CERCIGNANI:1988}. In the absence of sources and absorption it takes the form 
\begin{equation}
	\frac{1}{|v(\epsilon)|}\partial_t \psi(t,\epsilon,x,\Omega) + \Omega \cdot \nabla_x \psi(t,\epsilon,x,\Omega) = \mathcal{Q}(t,\epsilon,x) \left[ \psi(t,\epsilon,x,\Omega) \right]
	\label{eqn:linBoltzmann}
\end{equation}
where the unknown $\psi(t,\epsilon,x,\Omega)$ is the number density of particles, with respect to the measure $d\epsilon d\Omega dx$, located at $x$ and moving in direction $\Omega\in\US{2}$ with energy $\epsilon$ at time $t$. $|v(\epsilon)|$ is the absolute velocity of a particle with energy $\epsilon$. The scattering operator $\mathcal{Q}$ describes the change in time due to angular deflections and energy-loss as particles interact with the background medium through collisions.\\
Two situations that are of particular practical relevance allow to consider models of the same structure but reduced phase space.
\begin{enumerate}
	\item When the energy of particles is fixed we can omit the energy variable and obtain
	the \emph{radiative transfer equation}\cite{CASE:1967} (RT)
	\begin{equation}
		\partial_t \psi(t,x,\Omega) + \Omega \cdot \nabla_x \psi(t,x,\Omega) =  \underbrace{\int_{S^2} \sigma_s(x,\Omega'\cdot\Omega) \psi(t,x,\Omega')\, d\Omega' - \sigma_t(x) \psi(t,x,\Omega)}_{ \eqdef \mathcal{Q}^{RT}(x) \psi(t,x,\Omega)}
		\label{eqn:radiationTransport}
	\end{equation}
	where $\sigma_s$ is a non-negative scattering cross-section consisting of the density of scattering centers $N_V(x)$ in the background medium such that $\sigma_s(x,\cdot)/N_V(x)$ is the probability density of the angular deflection; $\sigma_t$ is the respective total cross-section. Note that $\sigma_s$, and thereby also $\sigma_t$, can depend on time.
	\item Many situations allow to neglect the time dependency and it can be assumed that particles loose a significant amount of energy only by a sequence of collisions with each collision changing the energy only slightly. In these situations one usually describes the particle system in terms of the particle fluence $\hat\psi(x,\epsilon,\Omega) \eqdef |v(\epsilon)| \psi(x,\epsilon,\Omega)$ and employs an evolution equation in energy space called \emph{Boltzmann equation in continuous slowing down approximation}\cite{LARSEN:1997} (BCSD)
	\begin{equation}
		-\partial_\epsilon \left( S(\epsilon,x) \hat\psi(\epsilon,x,\Omega)\right)  +\,  \Omega\cdot\nabla_x \,\hat\psi(\epsilon,x,\Omega) = \underbrace{\int_{S^2} \tilde\sigma_s(\epsilon,x,\Omega'\cdot\Omega) \,\hat\psi(\epsilon,x,\Omega')\,d\Omega'
		- \tilde\sigma_t(\epsilon,x) \,\hat\psi(\epsilon,x,\Omega)}_{ \eqdef \mathcal{Q}^{CSD}(\epsilon,x) \hat\psi(\epsilon,x,\Omega)} ,
		\label{eqn:BCSD}
	\end{equation}
	where the stopping power $S$ describes the average energy loss per distance traveled, $\tilde\sigma_{s}$ denotes a scattering cross-section involving elastic and inelastic collisions and $\tilde\sigma_t$ is the joint total scattering cross-section of elastic and inelastic collisions.\\
	Note that through the variable transformation $\epsilon(t) = \epsilon_{\text{max}} - t$ the BCSD can be transformed into a pseudo equation of radiation transport by setting the stopping power $S$ to one.
\end{enumerate}
From the viewpoint of this work, all three equations (\cref{eqn:linBoltzmann,eqn:radiationTransport,eqn:BCSD}) are of the same structure. We will focus on the equation for radiative transfer in the following analysis.\\

The \textbf{boundary conditions} for radiative transfer are of Dirichlet type and prescribe the incoming half of the particle distribution at the boundary:
\begin{equation}
	\psi(t,x,\Omega) \overset{!}{=}
		\psi_{in}(t,x,\Omega) \qquad \forall \Omega\in\US{2}: n\cdot\Omega < 0,
		\label{eqn:boltzmannBC}
\end{equation}
where $n$ denotes the outward-pointing normal vector at $x\in\partial G$ and, with $(\partial G \times {\US{2})}_{-} \eqdef \lbrace (x,\Omega)\in\partial G\times\US{2}: n\cdot\Omega < 0 \rbrace$, $\psi_{in}:\R_+\times(\partial G \times {\US{2})}_{-}\to\R_+$ is a given distribution of incoming particles. 
It is common sense, that the energy of particles in a bounded domain $G\subset\R^3$ at time $T$ is bounded from above by the energy of the particles that initially ($t=0$) occupied the domain augmented by the energy of the particles that entered through the boundary. This natural fact is reflected by the following \emph{energy bound} for a solution $\psi$ of of radiative transfer equation \eqref{eqn:radiationTransport}, see also \cite{Ringhofer:2001,Sarna:2018}
\begin{equation}
	\| \psi(T,\,\cdot\,) \|^2_{ L^{2}(G \times \US{2}) } 
	\leq 
	\| \psi_0 \|^2_{ L^{2}(G \times \US{2}) } 
	+ 
	\| \psi_{\text{in}} \|^2_{ L^{2}((0,T) \times (\partial G \times {\US{2})}_{-} )} .
	\label{eqn:energyStabilityRT}
\end{equation}
Note that this energy bound, due to linearity of \cref{eqn:linBoltzmann} in $\psi$, implies continuous dependency of solutions on the initial and boundary data given by $\psi_0$ and $\psi_{\text{in}}$, and assures uniqueness of solutions. For the particular case that no particles enter the domain $\psi_{\text{in}}(t,\,\cdot\,) \equiv 0$ at any time $t$, we can deduce that the number of particles cannot increase, i.e. the energy bound becomes
\begin{equation}
	\partial_t\, \| \psi(t,\,\cdot\,) \|^2_{ L^{2}(G \times \US{2})} \leq 0.
	\label{eqn:energyStabilityRTDecay}
\end{equation}
The energy bound \eqref{eqn:energyStabilityRT} reflects a fundamental property of the physics described by the radiative transfer equation and is relevant for the analysis of well posedness in terms of uniqueness and continuous dependency on data of solutions.
Therefore it is desirable to assure an analogous energy bound also for reduced models aiming to approximate the linear Boltzmann equation, such as the \PN{N} approximation considered in this work.\\

In this paper we present boundary conditions for the $P_{\text{N}}$ approximation of radiative transfer, which assure an energy bound analogous to \eqref{eqn:energyStabilityRT} for the approximate solution and are compatible with the \PN{N} equations, i.e. incoming characteristic waves are uniquely determined by outgoing characteristic waves and boundary data. We introduce the boundary conditions in two steps. In the first step (\cref{sec:analytic}), we investigate a generic first order system of linear transport equations with relaxation, an abstraction of \PN{N} equations, and show stability and compatibility for boundary conditions which have a particular structure. In the second step (\cref{sec:PNBC}), we derive the boundary conditions for $P_{\text{N}}$ equations and show that these are of the desired structure.
In \cref{sec:numericalSolution} we describe a discretization of $P_{\text{N}}$ equations in space using summation by parts finite difference on staggered grid for spatial derivatives and simultaneous approximation terms for the boundary conditions, show that a semi discrete solution satisfies an energy bound analogous to \eqref{eqn:energyStabilityRT} and present first numerical results in \cref{sec:numericalResults}.
	\section{Boundary conditions for Onsager Compatible Systems}
\label{sec:analytic}
In this section we investigate boundary conditions of a specific structure for an abstract form of transport equations.
\begin{problem}
	Let $\domain\subset\R^d$. Find solution $u:\mathbb{R}_+\times\domain \to \mathbb{R}^m$ of the first-order system
	\begin{alignat}{2}
		\partial_t u + \sum_{i=1}^d A^{(i)} \partial_{x_{i}} u &= P u \qquad && \forall x\in\domain,\ \forall\, t>0
		\label{eqn:genericPDE}
	\end{alignat}
	with transport matrices  $A^{(i)} \in \mathbb{R}^{m \times m}$ and negativ semi definit relaxation matrix  $P \in \mathbb{R}^{m \times m}$,
	for given initial datum $u_0:G\to\R^m$ 
	\begin{alignat}{2}
		u(0,\cdot) &= u_0 \qquad&& \forall x\in\domain \label{eqn:genericIC}.
	\end{alignat}
	\label{prob:GenericCauchyProblem}
\end{problem}
\noindent We augment Problem 1 with the following assumptions.
\begin{assumption}[Rotational Invariance]
	\Cref{eqn:genericPDE} is \emph{rotational invariant}, i.e. for a orthogonal matrix $\mathcal{R}\in\R^{d\times d}$ mapping a cartesian coordinate system $C$ onto $\hat C$ there is an orthogonal rotation matrix $Q(\mathcal{R}) \in \R^{m\times m}$ that commutes with the spatial operator $\mathcal{L} u \eqdef \sum_{i=1}^d A^{(i)} \partial_{x_{i}} u + P u$ in the sense that
	\begin{equation}
		\hat{\mathcal{L}} \circ Q(\mathcal{R}) = Q(\mathcal{R}) \circ \mathcal{L},
	\end{equation}
	where $\hat{\mathcal{L}} u \eqdef \sum_{i=1}^d A^{(i)} \partial_{\hat{x}_{i}} u + P u$ with $\partial_{\hat{x}_{i}}\, \cdot= \sum_{i=1}^d \mathcal{R}_{ji} \partial_{x_{j}}\, \cdot\, $.
	\label{Ass:RotationalInvariance}
\end{assumption}

The rotational invariance allows to express the directional transport matrix $A^{(n)} := \sum_{i=1}^d n_i A^{(i)}$ with a unit direction $n\in\US{d-1}$ in terms of a transport matrix $A^{(j)}$ through rotation
\begin{align}
	A^{(n)} &= Q(R_{n\to e_{j}}) A^{(j)} \left(Q(R_{n\to e_{j}})\right)^T
	\label{eqn:An1}
\end{align}
where $R_{n\to e_{j}}$ denotes a rotation that maps direction $n$ onto unit vector $e_j$. Vice versa, we assume $A^{(n)}$ to be given, use \cref{eqn:An1} to express $A^{(i)}$ in terms of $A^{(n)}$ and obtain
\begin{equation}
	\sum_{i=1}^d A^{(i)}\partial_{x_{i}}(\cdot)
	=\sum_{i=1}^d \left( Q(R_{e_{i}\to n}) A^{(n)} \left(Q(R_{e_{i}\to n})\right)^T \right) \partial_{x_{i}}(\cdot).
\end{equation}
Consequently, it suffices to characterize the transport operator by essentially one transport matrix $A^{(n)}$.

\begin{assumption}[Block Structure Symmetry/Onsager Compatibility]
	For at least one Cartesian direction $e_i$ the variables in $u$ can be ordered such that the transport matrix $A^{(i)}$ is of the form
	\begin{equation}
		A^{(i)} = \begin{pmatrix} 0 & \hat{A}^{(i)} \\ \left(\hat{A}^{(i)}\right)^T	& 0	\end{pmatrix}
		\label{eqn:OnsagerForm}
	\end{equation}
	where $\hat{A}^{(i)} \in\R^{r \times s}$ has full row rank, i.e. $r\leq s$ and $\rank(\hat A^{(i)})=r$.
\end{assumption}

We introduce compact notation for a splitting in $u$ with respect to a direction $n\in\mathbb{S}^{d-1}$:
\begin{align}
	\unp &\eqdef \left(\Qnp\right)^T u \qquad
	\unm \eqdef \left(\Qnm\right)^T u,
\end{align}
where $\Qnp\in\R^{m\times r}$ and $\Qnm\in\R^{m\times (m-r)}$ such that $Q(R_{n\to e_{j}})=[\Qnp,\Qnm]$ for some rotation $R_{n\to e_{j}}$.
The splitting is induced by the block structure of $A^{(n)}$ in the sense that
\begin{equation}
	u^T A^{(n)} u = \left(\unp\right)^T \hat{A} \unm + \left( \unm \right)^T \hat{A}^T \unp.
	\label{eqn:splittingByOnsagerStructure}
\end{equation}
Later we will refer to the variable sets $\unp$ and $\unm$ as odd and even. 
Note that the symmetry of the transport matrix $A^{(n)}$ also implies hyperbolicity of \cref{eqn:genericPDE}.

We define the class of boundary conditions investigated in this work, similar to the work done in \cite{Sarna:2018,Rana:2016,Struchtrup:2007}.
\begin{definition}[Onsager Boundary Condition]
	Let $n\in\mathbb{S}^{d-1}$ denotes the outward-pointing boundary normal at the respective boundary point $x\in\partial\domain$. We call a boundary condition \emph{Onsager boundary condition}, if it is of the form
	\begin{equation}
		\unp(x) = L \hat{A} \unm(x) + g(x)
		\qquad \forall x\in\partial\domain,
	\label{eqn:OBC}
	\end{equation}
	with 
	\begin{enumerate}[(i)]
		\item $L\in\mathbb{R}^{r\times r}$ symmetric semi positive definite and
		\item boundary source $g(x)\in\Ima(L)$.
	\end{enumerate}
	\label{def:onsagerBC}
\end{definition}

In the remainder of this section we show that OBCs imply an energy bound for a solution $u$ in terms of its data $u_0$ and $g$, and verify their compatibility with the hyperbolic characteristics of the field equation \eqref{eqn:genericPDE}. We start with proofing the energy bound, where we use $\| v \|_{L^2(D)}$ to denote the energy-norm of a vector valued function $v:D\to\R^s$ with $v_i\in L^2(D)$
\begin{equation}
	\| v \|_{L^2(D)}\eqdef \int_D \left(v(x)\right)^T v(x)\,dx.
\end{equation}

\begin{theorem}[Energy Bound]
	A solution $u$ to \cref{prob:GenericCauchyProblem} with Onsager boundary condition
	is bounded by the initial and boundary data
	\begin{align}
		\| u(T,\cdot) \|^2_{L^2(G)} \leq \| u_0 \|^2_{L^2(\domain)} 
		+ \| L^+ \|_2^2\,
		\| g \|^2_{L^2((0,T)\times\partial\domain},
		\label{eqn:energyStabilityOnsagerSystem}
	\end{align}
	where $\| L^+ \|_2$ denotes the operator norm of the pseudo inverse $L^+$ of $L$, i.e. $\| L^+ \|_2 = \max\lbrace \lambda^{-1}: \lambda\in\sigma(L)\setminus\lbrace 0 \rbrace \rbrace$.
	For the particular case that boundary source vanishes $g(T,\,\cdot\,)\equiv0$ at all times $t$ we get that energy cannot increase
	\begin{equation}
		\partial_t\, \| u \|_{L^2(G)} \leq 0.
	\end{equation}
	\label{thrm:energyStability}
\end{theorem}
\begin{proof}
	We multiply equation \eqref{eqn:genericPDE} with $u^T$ from the left, use the product rule of differentiation and integrate over the domain $G$ to obtain
	\begin{align*}
		\partial_t u + \sum_{i=1}^d A_i \partial_{x_{i}} u + P u &= 0
		\quad\Rightarrow\quad
		\partial_t \| u \|^2_{L^2(G)} = - \int_G \sum_{i=1}^d \partial_{x_{i}} ( u^T A_i u ) \,dx + 2 \int_G u^T P u \,dx.
	\end{align*}
	Applying the Gauss-theorem and using \cref{eqn:splittingByOnsagerStructure} we obtain
	\begin{align*}
		\partial_t \| u \|^2_{L^2(G)}
		&= -\int_{\partial G} u^T \underbrace{\left( \sum_{i=1}^d n_i A_i \right)}_{A^{(n)}} u \, dS + 2 \int_G u^T P u \,dx
		= - 2 \int_{\partial G} \left( \unp \right)^T \hat A \unm \, dS + 2 \int_G u^T P u \,dx.
	\end{align*}
	Now we use the Onsager boundary condition to express $\unp$ in terms of $\unm$ and $g$ to obtain
	\begin{align*}
		\partial_t \| u \|^2_{L^2(G)}
		&= - 2 \int_{\partial G} (L \hat{A} \unm + g)^T \hat A \unm \, dS + 2 \int_G u^T P u \,dx \\
		&= \int_{\partial G} \underbrace{-2 g^T \hat{A} \unm - (\hat{A} \unm)^T L (\hat{A} \unm)}_{\overset{(v\eqdef-\hat{A}\unm)}{=} 2g^Tv - v^T L v} \, dS 
		+ \int_{\partial G} \underbrace{ - (\hat{A} \unm)^T L (\hat{A} \unm)}_{\leq 0} \, dS + 2 \int_G \underbrace{ u^T P u }_{\leq 0} \,dx. 
	\end{align*}
	We now bound $2 g^T v - v^T L v$ from above by $g$ and $L$: As $g$ is in the image of $L$, we can always find a $v_g\in\R^c$ with $g=L v_g$ and with $L$ symmetric positive-definite we obtain $\forall v$
	\begin{align*}
		2g^T v - v^T L v &= {v_g}^T L v_g - (v-v_g)^T L (v-v_g) \leq |{v_g}^T L v_g| = |g^T v_g|  \\
		&\leq \| g \|_2 \| v_g \|_2 = \| g \|_2 \| L^+ g \|_2 
		\leq \| L^+ \|_2 \| g \|^2_2 = \| L^+ \|_2 \, g^T g  .
	\end{align*}
	With the above inequality
	\begin{align*}
		\partial_t \| u \|^2_{L^2(G)}
		&\leq \| L^+ \|_2 \int_{\partial G} g^T g \, dS 
		= \| L^+ \|_2 \| g(t,\,\cdot\,) \|^2_{L^2(\partial\domain)}
	\end{align*}
	and we obtain the energy estimate
	\begin{align*}
		\| u(T,\cdot) \|^2_{L^2(G)}& \leq \| u_0 \|^2_{L^2(G)} + 
		\| L^+ \|_2 \| \| g \|^2_{L^2((0,T) \times \partial \domain)}.
	\end{align*}
\end{proof}

It remains to show that OBCs are compatible with the characteristics of the field equations \eqref{eqn:genericPDE}, which we show using the following result on the eigenstructure of Onsager compatible matrices.

\begin{theorem}[Eigenstructure of Onsager Compatible Matrices]
	\label{thrm:eigenstructure}
	Let $A\in\mathbb{R}^{m\times m}$ be an Onsager compatible matrix such that
	\begin{equation*}
		A = \begin{pmatrix}
		0 & \hat A \\
		{\hat A}^T & 0
		\end{pmatrix}.
	\end{equation*}
	$\hat A\in\R^{r\times(m-r)}$. Then $A$ has an eigendecomposition $A=X \Lambda X^T$ with
	\begin{equation}
	\label{eqn:eigenstructure}
		X = \frac{1}{\sqrt{2}}
		\begin{pmatrix}
		\hat{X} & 0 & \hat{X} \\
		\tilde{X} & \sqrt{2} X^k & -\tilde{X}
		\end{pmatrix} \in\mathbb{R}^{m \times m}
		\qquad
		\Lambda =
		\begin{pmatrix}
		\Lambda_p & 0 & 0 \\
		0 & 0 & 0 \\
		0 & 0 & - \Lambda_p \\
		\end{pmatrix} \in\mathbb{R}^{m \times m}
	\end{equation}
	where
	\begin{enumerate}[(i)]
		\item $\Lambda_p \in \mathbb{R}^{r \times r}$ diagonal positive definit,
		\item $\hat{X} \in \mathbb{R}^{r \times r}$ orthogonal,
		\item $\tilde{X} \in \mathbb{R}^{(m-r) \times r}$ has orthonormal columns and
		\item the columns of $X^k \in \mathbb{R}^{(m-r) \times (m-2r)}$ form an orthonormal basis of $\ker(\hat{A})$ and are orthogonal to the columns of $\tilde{X}$.
	\end{enumerate}
\end{theorem}
\begin{proof}
	We can always write $\hat A$ by its singular value decomposition, say $\hat{A}=U\Sigma V^T$ where $U\in\mathbb{R}^{r \times r}, V\in\mathbb{R}^{(m-r) \times (m-r)}$ orthogonal and $\Sigma\in\mathbb{R}^{r \times (m-r)}$ with $\Sigma_{ij}=0$ for $i\neq j$. We identify $V=[V_1,V_2]$ with $V_1\in\mathbb{R}^{(m-r) \times r}$ and $V_2\in\mathbb{R}^{(m-r) \times ((m-2r)}$, and $\Sigma = [\Sigma_1, 0]$ with $\Sigma_1\in\mathbb{R}^{r \times r}$, such that $U\Sigma_1 V_1^T=\hat A$. Choosing $U=\hat{X},\tilde{X}=V_1,X^k=V_2,\Lambda_p=\Sigma_1$ yields the result.
\end{proof}

The structure of the diagonal matrix $\Lambda$ in \eqref{eqn:eigenstructure} tells us that the spectrum of the transport matrix $A^{(n)}$, is symmetric about the origin. We can deduce that at a boundary point we have $r$ incoming, $r$ outgoing and $(m-2r)$ standing characteristic waves. Existence and uniqueness of solutions to \cref{prob:GenericCauchyProblem} in the context of hyperbolic equations require that the boundary conditions uniquely determine the incoming, and only the incoming, characteristic waves. In order to allow an energy bound of the form in \eqref{eqn:energyStabilityRT} the incoming waves must be determined by the boundary data and outgoing waves, and be independent of incoming waves.
From the following theorem we can immediately deduce that Onsager boundary conditions satisfy the aforementioned requirements.

\begin{theorem}[Characteristic Compatibility]
	\label{thrm:charBC}
	Onsager boundary condition \eqref{eqn:OBC} is equivalent to the boundary condition
	\begin{equation}
		(L\hat{X}\Lambda_p + \hat{X}) w_- = (L\hat{X}\Lambda_p - \hat{X}) w_+ + \sqrt{2} {g}
		\label{eqn:chrBC}
	\end{equation}
	for the outgoing and incoming characteristic variables $w_+$ and $w_-$,  given by 
	\begin{equation}
		\begin{pmatrix}
			w_- \\
			w_+
		\end{pmatrix}
		=
		\frac{1}{\sqrt{2}}
		\begin{pmatrix}
			{\hat X}^T & {\tilde X}^T\\
			{\hat X}^T & -{\tilde X}^T
		\end{pmatrix}
		\begin{pmatrix}
		\unp \\
		\unm
		\end{pmatrix}.
	\end{equation}
	Matrix $(L\hat{X}\Lambda_p + \hat{X})$ is invertible and therefore Onsager boundary conditions uniquely determine the incoming characteristic waves by the outgoing waves and the boundary data.
\end{theorem}

\begin{proof}
	To prove the equivalence of the Onsager boundary condition \eqref{eqn:OBC} and the characteristic boundary condition \eqref{eqn:chrBC} we first express $\unp$ and $\unm$ by the characteristic variables 
	$$
		\begin{pmatrix}
			\unp \\ \unm
		\end{pmatrix}
		=
		\frac{1}{\sqrt{2}}
		\begin{pmatrix}
			\hat{X}  & 0 & \hat{X} \\
			\tilde{X} & \sqrt{2}X^k & -\tilde{X}
		\end{pmatrix}
		\begin{pmatrix}
			w_+ \\ w_0 \\ w_-
		\end{pmatrix}
		=
		\frac{1}{\sqrt{2}}
		\begin{pmatrix}
			\hat{X} (w_+ + w_-) \\
			\tilde{X} (w_+ - w_-) + \sqrt{2}X^k w_0
		\end{pmatrix}
	$$
	where $w_+$ denotes the outgoing, $w_-$ the incoming characteristics and $w_0$ the characteristics with zero velocity in normal direction. We replace $\unp$ and $\unm$ in the boundary condition by their expression in characteristic variables, factor out the characteristic variables and with $\hat{A}=\hat{X}\Lambda_p\tilde{X}^T$ we obtain
	\begin{align*}
		&\unp = L \hat{A} \unm + {g}\\
		\Leftrightarrow \quad 
		& \hat{X} (w_+ + w_-) = L (\hat{X}\Lambda_p\tilde{X}^T) ( \tilde{X} (w_+ - w_-) + \sqrt{2}X^k w_0) + \sqrt{2} g\\
		\Leftrightarrow \quad 
		& (L \hat{X}\Lambda_p + \hat{X}) w_- = (L \hat{X}\Lambda_p - \hat{X}) w_+ + \sqrt{2} g.
	\end{align*}
	The matrix $(L \hat{X}\Lambda_p + \hat{X})$ is invertible as we can write it as a product of invertible matrices
	$$
		(L \hat{X}\Lambda_p + \hat{X}) = \hat{X} \underbrace{(\hat{X}^T L \hat{X} + {\Lambda_p^{-1}})}_{\text{positiv definit}} \Lambda_p.
	$$
\end{proof}
 	\section{Stable Boundary Conditions for $\PN{N}$ equations}
\label{sec:PNBC}

In this section we present boundary conditions, which ensure that the \PN{N} approximation maintains the energy estimate \eqref{eqn:energyStabilityRT}. We will exploit the previous \cref{sec:analytic}: First, in \cref{ssec:PNapproximation}, we embed \PN{N} equations into the context of \cref{sec:analytic}, i.e. we show that \PN{N} equations are instances of the abstract field equations \eqref{eqn:genericPDE} in \cref{prob:GenericCauchyProblem}. Secondly, in \cref{ssec:PNBC}, we discuss the formulation of boundary conditions for \PN{N} equations and present Onsager boundary conditions.

\subsection{\PN{N} Approximation}
\label{ssec:PNapproximation}

\subsubsection{Motivation}
\label{ssec:MotivationPN}
Due to the high dimension of the number density $\psi$, in general $\psi$ is a function of six variables, solving radiation transport equation by a direct discretization in phase space is usually computationally too expensive. Besides the curse of dimensionality, the evolution of $\psi$ is non-local as the collision operator introduces a coupling in energy and direction space. 
Therefore it is usually necessary to use reduced models, one of which are so-called \PN{N} models, which we consider in this work. 

A \PN{N} model follows from a spectral Galerkin discretization in direction space using spherical harmonic functions as basis \cite{Davison:1958}. The number density $\psi$ is approximated by an expansion in spherical harmonic functions $Y_l^k:\US{2}\to\R$ up to order $N\in\N$
\begin{equation}
	\psi(t,x,\Omega) \approx \psi_{\PN{N}}(t,x,\Omega) 
	= \sum_{\substack{l\leq N, |k|\leq l}} u_l^k(t,x) Y_l^k(\Omega),
	\label{eq:PnApprox}
\end{equation}
where the expansion coefficients are the spherical harmonic moments of the angular distribution. Evolution equations for the coefficients $u_l^k$, the \PN{N} equations, are obtained by testing the equation of radiation transport with spherical harmonic functions up to order $N$
\begin{equation}
	\int_{S^2} (\partial_t \psi_{\PN{N}} + \Omega \cdot \nabla_x \psi_{\PN{N}}) Y_{l'}^{k'}\,d\Omega \
	= \int_{S^2} (\mathcal{Q}^{RT} \psi_{\PN{N}}) Y_{l'}^{k'}\,d\Omega
	\qquad
	\forall l'\leq N,\,\forall |k'|\leq l'.
	\label{eqn:PNeqnSketch}
\end{equation}
The phase space gets reduced by two dimensions at the cost of replacing a scalar equation by a system of equations. More importantly, the spherical harmonic functions are eigenfunctions of the scattering operator whenever the cross sections are isotropic, i.e. depend only on the deflection angle $(\Omega \cdot \Omega')$, which reduces the computational cost of the scattering operator in the Galerkin method, as the expansion coefficients decouple in direction space.

\subsubsection{Spherical Harmonic Functions and Notation}
\label{ssec:sphericalHarmonics}
Spherical harmonic functions form a \emph{complete orthonormal set} of the space $L^2(\US{2};\R)$. The real spherical harmonic $\SH{l}{k}:\US{2}\to\R$ of degree $l\in\N_0$ and order $k$ ($|k|\leq l$) is given by
\begin{equation}
	\SH{l}{k}(\Omega(\mu,\varphi))
	C_{l,|k|}
	P_l^{|k|}(\mu)
	\begin{cases}
		\cos(|k|\varphi)& k>0\\
		1/\sqrt{2}& k=0\\
		\sin(|k|\varphi)& k<0,
	\end{cases}
	\qquad\quad
	C_{l,|k|} = (-1)^{|k|}
 	\sqrt{\frac{2l+1}{2\pi}\frac{(l-|k|)!}{(l+|k|)!}}
\end{equation}
where $P_l^{|k|}$ denotes the associated Legendre function of respective degree and order, and $\mu$ and $\varphi$ are the polar and azimuthal angle of direction $\Omega$. These functions are well-studied, e.g. \cite{Atkinson:2014,Muller:1966}. We recapitulate properties essential to this work in a form that is customized to an observation we will use to classify spherical harmonics:
A spherical harmonic $\SH{l}{k}$ is either even or odd in the direction of a cartesian axes, where we call a function $f:S^2\to\mathbb{R}^n$ \emph{even/odd in direction $n\in \US{2}$}, if
\begin{equation}
	f\left(\Omega-2(n\cdot\Omega) n\right) = \raisebox{.2\height}{\scalebox{.8}{+}}
	\raisebox{.2\height}{\scalebox{.8}{/}}
	\raisebox{.14\height}{\scalebox{1.0}{-}}\,
	f(\Omega) \qquad \forall\, \Omega\in \US{2}.
	\label{eqn:ourEvenOddClass}
\end{equation}
Table \ref{tab:evenodd} gives simple formulas for the even/odd classification based on the degree $l$ and order $k$. 

When considering all three cartesian axes simultaneously the set of all spherical harmonics is separated into $2^3=8$ non-trivial sets, which can easily be verified using the formulas.
\begin{table}[t]
	\begin{center}
		\begin{tabular}{ccc}
			\toprule
			$n$ & $\SH{l}{k}$ even & $\SH{l}{k}$ odd \\
			\midrule
			$e_x$ 	& $\ \ (k<0 \wedge k\text{ odd} ) $ & $\ \ \ \ (k<0 \wedge k\text{ even} )$ \\
					& $\vee\ (k\geq 0 \wedge k \text{ even})$ & $\vee\ (k\geq 0 \wedge k \text{ odd})$ \\
			$e_y$ & $k \geq 0$ & $k<0$ \\ 
			$e_z$ & $(l+k)$ even & $(l+k)$ odd \\
			\bottomrule
		\end{tabular}
	\end{center}
	\caption{Even/Odd classification formulas for spherical harmonic functions with respect to the standard basis vectors $e_1$, $e_2$ and $e_3$ of the Cartesian coordinate system.}
	\label{tab:evenodd}
\end{table}
Note, that our classification is different from the usual even/odd classification of spherical harmonics based on the \emph{parity} referring to inversion about the origin, which solely depends on the even-odd property of the integer degree $l$, that is, $\SH{l}{k}(-\Omega) = \left(-1\right)^l\, \SH{l}{k}(\Omega)$.

For the sake of compact notation we introduce abbreviations for vector functions of spherical harmonics
\begin{alignat}{3}
& \SHnB{l}(\Omega) \eqdef (\SH{l}{-l}(\Omega),\SH{l}{-l+1}(\Omega),\ldots,\SH{l}{l}(\Omega))^T,
&\qquad&
\SHnB{l}:\US{2} \to \R^{2l+1}\\
& \SHn{N}(\Omega) \eqdef (\SHnB{0}(\Omega)^T,\SHnB{1}(\Omega)^T,\ldots,\SHnB{N}(\Omega)^T)^T.
&\qquad&
\SHN{N}:\US{2}\to\R^{(N+1)^2}
\end{alignat}
We use $\SHnB{l}^{e,i}(\Omega) \in \mathbb{R}^{l+1}$ to denote a vector function containing the spherical harmonic functions of degree $l$, that are even with respect to the direction $e_i$ ($i\in\lbrace 1,2,3\rbrace$) and denote the odd counterpart by $\SHnB{l}^{o,i}(\Omega) \in \mathbb{R}^{l}$. As above we also designate their collections up to degree $N$
\begin{alignat}{3}
& \SHn{N}^{e,i}(\Omega) \eqdef (\SHnB{0}^{e,i}(\Omega)^T,\ldots,\SHnB{N}^{e,i}(\Omega)^T)^T ,
&\qquad&
\SHn{N}^{e,i}:S^2 \to \R^{\frac{(N+2)(N+1)}{2}}\\
& \SHn{N}^{o,i}(\Omega) \eqdef (\SHnB{0}^{o,i}(\Omega)^T,\ldots,\SHnB{N}^{o,i}(\Omega)^T)^T ,
&\qquad& 
\SHn{N}^{o,i}:S^2 \to \R^{\frac{N(N+1)}{2}}.
\end{alignat}
As scalar producst for given $v\in L^2(\US{2};\mathbb{R}^{r})$ and $w\in L^2(\US{2};\mathbb{R}^{c})$ we use
$\langle v,w^T \rangle$ and $\langle v,w^T \rangle_{n_{\pm}}$ to denote the $(r\times c)$-matrices defined by
\begin{align*}
\left(\langle v,w^T\rangle\right)_{ij} \eqdef \int_{S^2} v_i\, w_j\, d\Omega 
\qquad
\left(\langle v,w^T\rangle_{n_{\pm}}\right)_{ij} \eqdef \int_{\Omega \cdot n \gtrless  0} v_i\, w_j\, d\Omega.
\end{align*}

We now recapitulate two properties which, alongside orthogonality, ensure that \PN{N} equations are instances of the abstract field equations \cref{eqn:genericPDE}. First, the function space spanned by spherical harmonics of degree $l$ is \emph{rotational invariant} in the sense that for a rotation $\mathcal{R}$ about the origin mapping a unit vector $\Omega$ onto $\Omega'=\mathcal{R}\Omega$ the following equality holds
\begin{equation}
	\SH{l}{k}(\mathcal{R} \Omega) = \SH{l}{k}(\Omega') = \sum_{k'=-l}^l D^l_{k'}(\mathcal{R}) \SH{l}{k'}(\Omega)
	\quad \text{with}\quad
	D^l_{k'}(\mathcal{R}) = \langle \SH{l}{k}(\mathcal{R}\,\cdot\,), \SH{l}{k'}(\,\cdot\,) \rangle.
	\label{eqn:rotationalInvarianceSH}
\end{equation}
We can immediately deduce that for the vector of spherical harmonics up to order $N$
\begin{equation}
	\SHN{N}(\mathcal{R} \Omega) = D(\mathcal{R}) \SHN{N}(\Omega) \quad\text{with}\quad D(\mathcal{R})\eqdef \langle \SHN{N}(\mathcal{R}\, \cdot\,), \SHN{N}(\,\cdot\,)^T \rangle.
	\label{eqn:rotationMatrixSH}
\end{equation}
Secondly, spherical harmonics fulfill \emph{recursion relations}, which, in order to emphasize the coupling between even and odd functions, we write as
\begin{subequations}
\begin{align}
	\Omega_i \SHnB{l}{o,i}(\Omega) &= \hat{A}^{(i)}_{l,-} \SHnB{l-1}{e,i}(\Omega) + \hat{A}^{(i)}_{l,+} \SHnB{l+1}{e,i}(\Omega)
	\label{eqn:recursion1}\\
	\Omega_i \SHnB{l}{e,i}(\Omega) &= (\hat{A}^{(i)}_{l-1,+})^T \SHnB{l-1}{o,i}(\Omega) + (\hat{A}^{(i)}_{l+1,-})^T \SHnB{l+1}{o,i}(\Omega) \label{eqn:recursion2}
\end{align}
\end{subequations}
with
\begin{equation}
	\hat{A}^{(i)}_{l,\pm} = \langle \Omega_i \SHnB{l}{o,i}, (\SHnB{l$\pm$1}{e,i})^T \rangle \in \mathbb{R}^{l\times(l\pm 1 + 1)}.
	\label{eqn:recusionMatrix}
\end{equation}
In order to point out the specific structure of the recursion relation more clearly, we schematically write out the recursion relation for all spherical harmonics. Gathering all even and odd spherical harmonics, we can write the recursion relation as
\begin{equation}
	\Omega_i
	\begin{pmatrix}
		\SHno{$\infty$}{i}(\Omega)\\
		\SHne{$\infty$}{i}(\Omega)
	\end{pmatrix}
	=
	\begin{pmatrix}
	0 & \hat{A}^{(i)} \\
	\left(\hat{A}^{(i)}\right)^T & 0
	\end{pmatrix}
	\begin{pmatrix}
		\SHno{$\infty$}{i}(\Omega)\\
		\SHne{$\infty$}{i}(\Omega)
	\end{pmatrix}
	\label{eqn:recursionRelationCoupling}
\end{equation}
with the infinite matrix
\begin{equation}
	\hat{A}^{(i)}
	= 
	\langle \Omega_i \SHno{$\infty$}{i}, \left(\SHne{$\infty$}{i}\right)^T \rangle
	=
	\begin{pmatrix}
		\hat{A}^{(i)}_{1,-}&0&\hat{A}^{(i)}_{1,+}&0&\dots\\
		0&\hat{A}^{(i)}_{2,-}&0&\hat{A}^{(i)}_{2,+}&\ddots\\
		\vdots&\ddots&\ddots&\ddots&\ddots
	\end{pmatrix}.
	\label{eqn:recursionRelationAhat}
\end{equation}

\subsubsection{\PN{N} Equations}
\label{ssec:PN}
Following (\ref{eq:PnApprox}) we replace $\psi$ in the equations of radiation transport \eqref{eqn:radiationTransport} by an expansion in spherical harmonics up to order $N$
$$
	\psi_{\PN{N}}(t,x,\Omega) 
	= u^T(t,x) \SHN{N}(\Omega),
$$
test the equation \eqref{eqn:radiationTransport} by multiplication with spherical harmonics up to order $N$, and integrate over the unit sphere 
\begin{align}
	\langle \partial_t (u^T \SHN{N}), \SHN{N} \rangle 
	+ \langle \Omega\cdot\nabla_x (u^T \SHN{N}), \SHN{N} \rangle
	&=  \langle \mathcal{Q}_{RT} (u^T \SHN{N}), \SHN{N} \rangle.
	\label{eqn:PNintermediate}
\end{align}
\Cref{eqn:PNintermediate} is equivalent to
\begin{align}
	\underbrace{
		\langle \SHN{N},\SHN{N}^T \rangle
	}_{=I_N} 
	\partial_t u  +
	\sum_{i=1}^3 
	\underbrace{
		\langle \Omega_i \SHN{N}, \SHN{N}^T \rangle
	}_{=:A^{(i)}_N} 
	\partial_{x_{i}} u 
	&= 
	\underbrace{
		\langle \mathcal{Q}_{RT} \SHN{N}, \SHN{N}^T \rangle
	}_{=:Q^{RT}_N} 
	u,
\end{align}
with $(\mathcal{Q}_{RT} \SHN{N})_i=\mathcal{Q}_{RT} (\SHN{N})_i$. Using the  orthonormality and recursion relation of spherical harmonics, as indicated by underbraces, we can confirm that \PN{N} equations can be written in the abstract form of \cref{eqn:genericPDE}.

\emph{Rotational invariance} of the \PN{N} equations follows from the rotational invariance of spherical harmonics, \cref{eqn:rotationalInvarianceSH}. It is obvious that the scattering matrix $Q_N$ and the advection matrices $A_N^{(i)}$ do not depend on the orientation of the coordinate system. Now let $\mathcal{R}\in\R^3$ be the matrix rotating a cartesian coordinate system and $D_N(\mathcal{R})$ the respective rotation matrix as defined in \cref{eqn:rotationMatrixSH}. It is easy to verify that
\begin{equation}
	Q_N D_N(\mathcal{R}) = D_N(\mathcal{R}) Q_N
	\label{eqn:commutationScattering}
\end{equation}
using the invariance of integrals over the unit sphere under rotation, orthogonality of the rotation matrices $\mathcal{R}$ and $D_N(\mathcal{R})$ and that scattering cross sections depend on the pre- and post-collision direction solely via the deflection angle. Likewise the transport operator commutes with the rotation $D_N(\mathcal{R})$
\begin{equation}
	\sum_{i=1}^d A_N^{(i)} \partial_{\tilde{x}_{i}} (D(\mathcal{R})\ \cdot\ ) 
	=
	D(\mathcal{R}) \sum_{i=1}^d A_N^{(i)} \partial_{x_{i}}(\ \cdot\ ),
	\label{eqn:commutationTransport}
\end{equation}
which can be verified using the same arguments and $\partial_{\tilde{x}_{i}}(\,\cdot\,) = \sum_{s=1}^d \mathcal{R}_{is} \partial_{x_{s}}(\,\cdot\,)$.\\

\emph{Block Structure Symmetry / Onsager compatibility} of the \PN{N} equations immediately follows from the structure of the recursion relations \cref{eqn:recursion1,eqn:recursion2} regarding the even/odd classification of spherical harmonics as in \cref{eqn:ourEvenOddClass}. As indicated by \cref{eqn:recursionRelationCoupling} it is clear that any permutation gathering first odd and then even components of $\SHN{N}$ allows to transform $A^{(i)}_N$ into the form of \eqref{eqn:OnsagerForm}. The full row rank requirement is satisfied, as the elements of $\Omega_i \SHNo{N}{i}(\Omega)$ are linearly independent.

\subsection{Boundary Conditions}
\label{ssec:PNBC}
While we know how to formulate boundary conditions for the linear kinetic equation, the right formulation of boundary conditions for \PN{N} equations is not clear. We cannot expect $\psi_{\PN{N}}$ to fulfill the original kinetic boundary conditions \cref{eqn:boltzmannBC} for general $\psi_{in}$, that is $\psi_{\PN{N}}$ can not be of the form
\begin{equation}
	\psi_{\PN{N}}(\Omega) =
	\begin{cases}
	\psi_{\PN{N}}(\Omega) & n\cdot\Omega \geq 0 \\
	\psi_{in}(\Omega) & n\cdot\Omega< 0,
	\end{cases}
	\label{eqn:boltzmannBCforPN}
\end{equation}
at the boundary, because the right hand side is discontinuous in general, while $\psi_{\PN{N}}$ is continuous. From \cref{thrm:eigenstructure} we know that the boundary conditions have to determine $N(N+1)/2$ incoming characteristic waves and hence, the right number of independent boundary conditions is $N(N+1)/2$. 

Two different strategies to transfer the kinetic boundary condition \eqref{eqn:boltzmannBC} into boundary conditions for \PN{N} equations are being used in the literature: 
Enforcing that the \PN{N} approximation $\psi_{\PN{N}}$ matches $\psi_{\text{in}}$ on a set of inward-pointing directions $\Omega_i$ ($\Omega_i \cdot n < 0$)
$$
	\psi_{\PN{N}}(\Omega_i) = \psi_\text{in}(\Omega_i) \qquad i=1,\ldots, N(N+1)/2,
$$
and secondly, equating ingoing half moments of the known distribution $\psi_{\text{in}}$ and the \PN{N} approximation $\psi_{\PN{N}}$ for a certain set of test function $\phi_i$
$$
	\int_{\Omega \cdot n<0} \psi_{\PN{N}}\, \phi_i\,d\Omega = \int_{\Omega \cdot n<0} \psi_\text{in}\, \phi_i\,d\Omega \qquad i=1,\ldots,N(N+1)/2.
$$
The two types of boundary conditions are referred to as Mark\cite{Mark:1944a,Mark:1944b} and Marshak\cite{Marshak:1947} type boundary conditions. In this work we will consider only boundary conditions of Marshak type with test functions $\phi_i\in\spanVecSpace{\SH{l}{k}:|k|\leq l\leq N}$. Rather than the strategy, it is the set of test functions or directions which is crucial, as the set determines the boundary conditions for the expansion coefficients.\\
From our previous observations regarding well-posedness, it is clear that the choice of test functions is crucial from a mathematical perspective. However, it is important to not forget, that the primary objective of \PN{N} equations and the accompanying boundary conditions is to serve as an accurate physical model and hence, the choice of test functions should also be reasonable from a modeling point of view. This is however difficult to check, as, apart from the macroscopic density $u_0^0$ and fluxes $u_1^i$, we lack a macroscopic interpretation to a general harmonic moment $u_l^k$.
We briefly discuss classical choices of test functions.

\subsubsection{Classical Marshak Boundary Conditions}
\label{ssec:PNBCclassic}
The common choice of test functions, which results in the right number of boundary conditions, is to test with the spherical harmonics of even (odd) parity, if the approximation order $N$ is odd (even). These choices are attractive for arbitrarily oriented boundary normal $n$, as the parity is invariant under rotation and the test space thereby independent of $n$. This approach of choosing the test functions is nevertheless questionable. For one, completely exchanging the test space when in- or decreasing the order $N$ by one is somehow unnatural, considering that \PN{N} equations follow from a truncated series expansion and therefore from a hierarchical family of models. And secondly, questions arise from the viewpoint of a physical model. The total incoming particle flux is probably the most fundamental quantity we want to capture exactly within a boundary model. This is not assured, when testing with a finite number of even parity harmonics. When testing with spherical harmonics of odd parity, the correct total incoming particle flux is enforced indirectly by the test functions $\SH{1}{-1}$, $\SH{1}{0}$ and $\SH{1}{1}$. In \cite{Williams:1971} this observation for the case of vacuum boundary condition is used to justify testing with odd parity spherical harmonics, i.e. $\phi_i\in\lbrace\SH{l}{k}:l \text{ odd}\rbrace$. This is not a valid argument for choosing all spherical harmonics of odd parity, as the correct total incoming particle flux can be imposed by a single equation (choose $\phi_i=\Omega\cdot n$).\\
To our knowledge, the issue of stability is taken into account only in the works \cite{Egger:2012} and \cite{Schlottbom:2015}, where stable variational formulations for \PN{N} equations are derived. These works however do not provide an explicit form of boundary conditions for the strong formulation of \PN{N} equations and rely on the parity splitting of spherical harmonics.

\subsubsection{Onsager Boundary Conditions}
\label{ssec:PNBConsager}
In the following we present a modification of Marshak type boundary conditions, which can be seen as a truncation of the recursion relation, consistent with the truncation of the series expansion, and show that the modification leads to Onsager boundary condition. We first motivate our choice of test functions in the Marshak conditions and translate the integral formulation into explicit conditions for the expansion coefficients.

We consider boundary conditions that are based on enforcing the original Boltzmann boundary condition \eqref{eqn:boltzmannBCforPN} weakly on the subspace of functions that are odd with respect to the boundary normal $n$, which we write as
\begin{equation}
\langle \psi_{\PN{N}}, \SHNo{N}{n} \rangle =
\langle \psi_{\PN{N}}, \SHNo{N}{n} \rangle_{n_{+}} + \langle \psi_{in}, \SHNo{N}{n} \rangle_{n_{-}}.
\label{eqn:BCPNabstract}
\end{equation}
Note that \Cref{eqn:BCPNabstract} is equivalent to the classical formulation of Marshak type boundary conditions, as the case $n\cdot\Omega\geq0$ in \cref{eqn:boltzmannBCforPN} is superfluous and hence the integral equation is governed by integrals over the incoming half sphere.
This test space is reasonable from various points of views, and compared to the previously described classical choice more natural:
\begin{itemize}
	\item the number of test functions is correct for all approximation orders $N$, in particular independent on whether $N$ is even or odd,
	\item the boundary conditions are hierarchical in $N$, that is, they form a cascade analogously to the hierarchy of $\PN{N}$-equations.
	\item in the limit $N\to\infty$ the test space forms a complete set on the ingoing half sphere $\lbrace \Omega\in\US{2}: \Omega\cdot n < 0 \rbrace$, hence the boundary condition is satisfied weakly in the limit,
	\item the correct total incoming particle flux, as most fundamental model property, is incorporated explicitly by the first equation and
	\item they allow to obtain Onsager boundary conditions and hence the desired consistency with the \PN{N} equations and an energy bound for solutions.
\end{itemize}

First, we rewrite \cref{eqn:BCPNabstract} as an equation for the expansion coefficients by splitting $\psi_{\PN{N}}$ into even and odd part $\psi_{\PN{N}}={\SHNo{N}{n}}^T  u^o + {\SHNe{N}{n}}^T u^e$, using the orthogonality property and replacing half sphere integrals of even functions by integrals over the full sphere
\begin{align}
	&\langle {\SHNo{N}{n}}^T  u^o + {\SHNe{N}{n}}^T u^e, \SHNo{N}{n} \rangle = 
	\langle {\SHNo{N}{n}}^T  u^o + {\SHNe{N}{n}}^T u^e, \SHNo{N}{n} \rangle_{n_{+}} + \langle \psi_{in}, \SHNo{N}{n} \rangle_{n_{-}}\\
	\Leftrightarrow \quad
	&\langle {\SHNo{N}{n}}, {\SHNo{N}{n}}^T \rangle u^o = 
	\frac{1}{2}\langle {\SHNo{N}{n}}, {\SHNo{N}{n}}^T \rangle u^o
	+ \langle {\SHNo{N}{n}}, {\SHNe{N}{n}}^T \rangle_{n_{+}} u^e 
	+ \langle \psi_{in}, \SHNo{N}{n} \rangle_{n_{-}}\\
	\Leftrightarrow \quad
	&\underbrace{\langle {\SHNo{N}{n}}, {\SHNo{N}{n}}^T \rangle}_{I} u^o = 
	\underbrace{2 \langle {\SHNo{N}{n}}, {\SHNe{N}{n}}^T \rangle_{n_{+}}}_{=:\tilde M} u^e 
	+ \underbrace{2 \langle \psi_{in}, \SHNo{N}{n} \rangle_{n_{-}}}_{=:g}\\
	\Leftrightarrow \quad
	&u^o = \tilde M u^e + g.
	\label{eqn:BCwithMcomplete}
\end{align}
We will point out a specific structure in the matrix $\tilde M$, which allows us to apply a natural modification to $\tilde M$ leading to Onsager boundary conditions. To keep the notation simple we consider the case $n=e_i$, i.e. the boundary normal $n$ coincides with the cartesian basis vectors $e_i$. The general case follows from rotation. We first insert an artificial $1=\frac{\Omega\cdot n}{\Omega\cdot n}=\frac{\Omega_i}{\Omega_i}$ for fixed $i$ of which we assign the nominator to the odd part and the denominator to the even
\begin{align}
	\tilde M = 2 \int_{\Omega_i \geq 0}\SHNo{N}{i}(\Omega) \left(\SHNe{N}{i}(\Omega)\right)^T\, d\Omega = 2 \int_{\Omega_i \geq 0}\frac{1}{\Omega_i}\SHNo{N}{i}(\Omega) \left( \Omega_i \SHNe{N}{i}(\Omega)\right)^T\, d\Omega.
	\label{eqn:Mcomplete}
\end{align}
Note that the term $\frac{1}{\Omega_i}\SHNo{N}{i}(\Omega)$ is finite for $\Omega_i \to 0$ as $\SHNo{N}{i}(\Omega)$ is at least first order in $\Omega_i$. The term $\Omega_i \SHNe{N}{i}(\Omega)$ is odd and can be written in terms of $\SHNo{N+1}{i}$ using the recursion relation \eqref{eqn:recursion2}
\begin{align}
	\Omega_i \SHNe{N}{i}(\Omega)
	&= \langle \Omega_i \SHNe{N}{i}(\Omega), \left(\SHNo{N+1}{i}(\Omega)\right)^T \rangle \SHNo{N+1}{i}(\Omega)
	= \left( \hat{A}^{(i)} \right)^T
	\SHNo{N}{i}(\Omega)
	+ R_{l=N+1} \SHnB{l=N+1}^{\text{o,i}}(\Omega)
\end{align}
The last term describes the contribution of higher order moments ($l>N$) in the recursion relation. Similar to the truncation of the series expansion, we truncate the recursion relation by neglecting the contribution of higher order moments to obtain an approximation of the recursion relation involving only terms up to order $N$
\begin{align}
	\Omega_i \SHNe{N}{i}(\Omega) \approx \left( \hat{A}^{(i)} \right)^T \SHNo{N}{i}(\Omega).
	\label{eqn:ApproxRecRel}
\end{align}
Inserting this approximation into \cref{eqn:Mcomplete} we obtain
\begin{equation}
	\tilde M \approx M = 2 \int_{\Omega_i \geq 0}\frac{1}{\Omega_i}\SHNo{N}{i}(\Omega) \left( \left( \hat{A}_N^{(i)} \right)^T \SHNo{N}{i}(\Omega) \right)^T\, d\Omega
	=
	\underbrace{2 \int_{\Omega_i \geq 0}\frac{1}{\Omega_i}\SHNo{N}{i}(\Omega) \left( \SHNo{N}{i}(\Omega) \right)^T\, d\Omega}_{=: L_N^{(i)}}\, \hat{A}_N^{(i)},
	\label{eqn:Mapprox}
\end{equation}
where we observe that the resulting matrix has the form of Onsager boundary condition \eqref{eqn:OBC}, i.e. the matrix is written as a symmetric positive-definite matrix $L_N$ multiplied by the flux matrix $\hat{A}_N^{(i)}$ and hence we obtain from (\ref{eqn:BCwithMcomplete}) the desired form of boundary conditions
\begin{equation}
	u^{o,i} = L_N^{(i)} \hat{A}_N^{(i)}\, u^{e,i} + g.
	\label{eqn:PNBConsagerForm}
\end{equation}
Note that $\tilde M$ in \cref{eqn:BCwithMcomplete} differs from $M$ only in the columns that correspond to even expansion coefficients of highest order ($l=N$). The effect of replacing $\tilde M$ by $M$ is therefore minor whenever the highest order moments are small, i.e. when the \PN{N} approximation, being based on a truncated series expansion, is a reasonable approximation in the first place. Note also, that $L_N^{(i)}$ is invertible and hence the condition $g(x)\in\Ima\left(L^{(i)}_N\right)$ is always satisfied. Therefore our boundary condition ensures the energy bound \eqref{eqn:energyStabilityOnsagerSystem} for general boundary data. 

 	\section{Stable Discretization of $P_N$ Equation}
\label{sec:numericalSolution}

In this section we describe the summation by parts (SBP)\cite{Kreiss:1974,Strand:1994,Olsson:1995,David:2014} finite difference discretization (FD) in space for \PN{N} equations with simultaneous approximation terms (SAT)\cite{Mark:1994} for Onsager boundary conditions.
Similar to the discretization of two dimensional \PN{N} equations with periodic or constant extrapolation boundaries conditions presented in \cite{Seibold:2014}, we exploit the coupling structure 
to discretize moments on staggered grids.
The discretization allows to assure an discrete analogue to the energy bound in \cref{thrm:energyStability} for general Onsager boundary conditions.
We first describe the discretization in space using the summation by parts (SBP) finite difference method on staggered grids, then describe the inclusion of the boundary condition using simultaneous approximation terms (SAT) in a manner that ensures stability of the semi-discrete system. We will neglect the relaxation term as its discretization is straight forward and, analogous to the continuous system, the conditions assuring energy stability of the semi-discrete system are independent of the relaxation term.

\subsection{Energy Stable Discretization in 1D}
We first describe the discretization in space and implementation of boundary conditions for the following one-dimensional problem on the an interval $I=\left( x_L, x_R\right)\subset\R$.
\begin{problem}
	\label{prob:1dAdvection}
	Find $u^o:\R_+\times I\to\R^{n_o}$ and $u^e:\R_+\times I\to\R^{n_e}$ such that 
	\begin{alignat}{2}
		\partial_t 
		\begin{pmatrix}
			u^o \\ u^e
		\end{pmatrix}
		+
		\begin{pmatrix}
		0 & \hat{A} \\ \hat{A}^T & 0
		\end{pmatrix}
		\partial_x
		\begin{pmatrix}
		u^o \\ u^e
		\end{pmatrix}
		&=0
		\qquad
		&\forall x\in I,\ \forall t>0
		\label{eqn:1DTransportSystem}
	\end{alignat}
	with Onsager boundary conditions
	\begin{alignat}{2}
		u^o &= -L \hat{A}\, u^e + g \qquad &x=x_L\\
		u^o &= L \hat{A}\,  u^e + g &x=x_R
	\end{alignat}
	and a compatible initial condition $u^o(0,\,\cdot\,)=u^o_0,\ u^e(0,\,\cdot\,)=u^e_0$.
\end{problem}

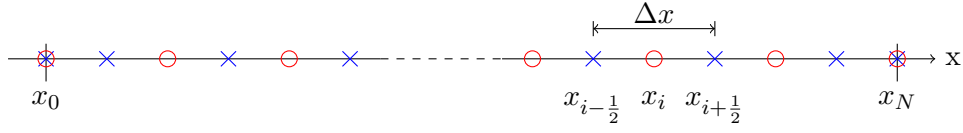
\begin{figure}
	\centering
	\begin{tikzpicture}
	\def\h{1.6}
	\def\a{0}
	\def\no{7}
	\def\ne{6}
	\def\delta{0.1}
	\def\b{\no*\h}
	\draw[-] (\a-0.5,0) -- (\a+2.75*\h,0);
	\draw[dashed] (\a+2.75*\h,0) -- (\a+3.75*\h,0);
	\draw[->] (\a+3.75*\h,0) -- (\b+0.5,0) node[right] {x};
	\draw[-] (\a,0.2) -- (\a,-0.3);
	\draw[-] (\b,0.2) -- (\b,-0.3);
	\foreach \i in {0,...,2}
	{
		\draw[color=blue] (\i*\h+0.5*\h-\delta,-\delta) -- (\i*\h+0.5*\h+\delta,+\delta);
		\draw[color=blue] (\i*\h+0.5*\h-\delta,+\delta) -- (\i*\h+0.5*\h+\delta,-\delta);
	}
	\foreach \i in {4,...,\ne}
	{
		\draw[color=blue] (\i*\h+0.5*\h-\delta,-\delta) -- (\i*\h+0.5*\h+\delta,+\delta);
		\draw[color=blue] (\i*\h+0.5*\h-\delta,+\delta) -- (\i*\h+0.5*\h+\delta,-\delta);
	}
	\foreach \i in {0,\no}
	{
		\draw[color=blue] (\i*\h-\delta,-\delta) -- (\i*\h+\delta,+\delta);
		\draw[color=blue] (\i*\h-\delta,+\delta) -- (\i*\h+\delta,-\delta);
	}
	\foreach \i in {0,...,2}{	\draw[color=red] (\i*\h,0) circle (3pt); }
	\foreach \i in {4,...,\no}{	\draw[color=red] (\i*\h,0) circle (3pt); }
	\node at (0,-0.55) {$x_0$};
	\draw[-] (4.5*\h,0.5) -- (4.5*\h,0.3); 
	\draw[-] (5.5*\h,0.5) -- (5.5*\h,0.3);
	\draw[<->] (4.5*\h,0.4) -- (5.5*\h,0.4);
	\node at (5.*\h,0.6) {$\Delta x$};
	\node at (4.5*\h,-0.65) {$x_{i-\frac{1}{2}}$};
	\node at (5*\h,-0.55) {$x_i$};
	\node at (5.5*\h,-0.65) {$x_{i+\frac{1}{2}}$};
	\node at (\no*\h,-0.55) {$x_N$};
\end{tikzpicture}
	\caption{Illustration of one-dimensional staggered grid: Red circles accommodate even variables (even grid) and blue crosses accommodate odd variabes (odd grid). At the boundary nodes both variables are present.}
	\label{fig:staggeredGrid1D}
\end{figure}

We discretize the interval $I$ by two grids, which are staggered to one another in the interior, as shown in Figure \ref{fig:staggeredGrid1D}.
We refer to the two grids as even and odd grid; identify the vectors of odd and even grid nodes by $x^o=(x_0,x_1,\ldots,x_N)^T\in\R^{N+1}$ and $x^e=(x_0,x_\frac{1}{2},x_{1+\frac{1}{2}},\ldots,x_{N-\frac{1}{2}},x_N)^T\in\R^{N+2}$; and denote the sets of grid nodes by $I_h^o\subset I$ and $I_h^e\subset I$. We use $f_{h_{\star}}$ to denote the grid function approximating the restricting of $f:I\to \R^n$ to $I_h^\star$, identify function values by subscripts $f_{h_{\star},i} \eqdef f_{h_{\star}}(x_i)$ and abuse $f_{h_{\star}}$ to also denote the vector obtained by concatenating the function values of all nodes, e.g.
\begin{equation}
	f_{h_{\star}} \eqdef \left({f_{h_{\star},0}^o}^T,\ldots,{f_{h_{\star},N}}^T\right)^T.
\end{equation}
The common notation 
\begin{equation}
	\langle f_h, g_h \rangle = \sum_{i} \left( f_{h,i} \right)^T g_{h,i}
\end{equation}
is used for the scalar product.

For the sake of a compact notation we introduce notation for operations on grid functions that can be interpretated as the product of a matrix with the vector of node values: Given two grids consisting of the nodes $I_h^1$ and $I_h^2$, with $|I_h^1|=N_1+1$ and $|I_h^2|=N_2+1$, and a matrix $M\in\R^{(N_2+1)\times (N_1+1)}$ we use $\mathcal{M}$ to denote the mapping of grid functions on $I_h^1$ to grid functions on $I_h^2$ given by
\begin{align}
	(\mathcal{M} f_{h_{1}})_i &\eqdef \sum_{j=0}^{N_1} M_{ij} f_{h_{1},j} \qquad i=0,\ldots,N_2
\end{align}
Furthermore we use $\mathcal{M}\mathcal{N}$, $\mathcal{M}\pm\mathcal{N}$, $\mathcal{M}^{-1}$ and $\mathcal{M}^T$ to denote the mappings defined by the matrices $MN$, $M+N$, $M^{-1}$ and $M^T$ whenever these are well defined.\\

The finite difference (FD) operators on staggered grids, which we consider here, have the character of a projection, in the sense that they approximates the first order spatial derivative of an odd/even grid function on the even/odd grid. In particular, we consider the special class of summation by parts (SBP) FD operators, which are constructed such that integration by parts can be imitated on the discrete level \cite{Oreilly:2017}. 
\begin{definition}[Staggered SBP-FD]
	 We call $\mathcal{D}^{e\to o}$ and $\mathcal{D}^{o\to e}$ SBP FD operators on the staggered grids $x^o$ and $x^e$, if they can be identified by matrices of the form
	\begin{alignat*}{2}
		D^o\eqdef D^{e\to o} &= \left(P^{o}\right)^{-1} Q^{o}, \qquad
		D^e\eqdef D^{o\to e} &= \left(P^{e}\right)^{-1} Q^{e}, 
	\end{alignat*}
	with $P^{o}\in\R^{(N+1)\times(N+1)}$, $P^{e}\in\R^{(N+2)\times(N+2)}$, $Q^{o}\in\R^{(N+1)\times(N+2)}$ and $Q^{o}\in\R^{(N+2)\times(N+1)}$, such that
	\begin{enumerate}[(i)]
		\item the order of \emph{accuracy} is $s\in\N$ for boundary nodes and 2$s$ for interior nodes, e.g.
		\begin{align*}
			\left( \mathcal{D}^{e\to o} (x^n)_{h_{e}} \right)_i &= \left( n x^{n-1} \right)_{h_{o},i}
			\qquad
			\begin{aligned}
				n&=0,\ldots,s  && \text{if $i$ is boundary node}\\
				n&=0,\ldots,2s && \text{if $i$ is internal node},
			\end{aligned}
		\end{align*}
		\item the so-called \emph{SBP-property}
		\begin{equation*}
			Q^o + \left(Q^e\right)^T = B^o = \left(B^e\right)^T = B \qquad\mbox{with}\qquad \langle f_{h_{o}}, \mathcal{B} g_{h_{e}} \rangle = {f_{h_{o},N}}^T g_{h_{e},N} - {f_{h_{o},0}}^T g_{h_{e},0}
		\end{equation*}
		is fullfilled and
		\item the matrices $P^{o}$ and $P^{e}$ are diagonal positive definite and thereby introduce \emph{discrete norms}
		\begin{align*}
			\| f_{h_{o}} \|^2_{h_{o}} \eqdef \langle f_{h_{o}}, \mathcal{P}^o f_{h_{o}} \rangle, \qquad 
			\| f_{h_{e}} \|^2_{h_{e}} \eqdef \langle f_{h_{e}}, \mathcal{P}^e f_{h_{e}} \rangle.
		\end{align*}
	\end{enumerate}
\end{definition}
For a detailed description of the construction of the matrices $P^{o}$, $P^{e}$, $Q^{o}$ and $Q^{e}$ we refer the reader to \cite{Oreilly:2017}.\\

We now discretize system \eqref{eqn:1DTransportSystem} in space by restricting  $u^o$ to $I_h^o$ and $u^e$ to $I_h^e$ and replacing the spatial derivatives of $(\hat{A} u^e)$ and $(\hat{A}^T u^o)$ by SBP-FD operators
\begin{nalign}
	\label{eqn:semiDiscreteSystem01}
	\partial_t u_{h_{o}}^o &+ \mathcal{D}^{e\to o} (\hat{A} u^e)_{h_{e}}=0 \\
	\partial_t u_{h_{e}}^e &+ \mathcal{D}^{o\to e} (\hat{A}^T u^o)_{h_{o}}=0.
\end{nalign}
SBP-FD operators allow to mimic integration by parts of the continuous system, as
\begin{align*}
	\frac{1}{2} \partial_t \| u_{h_{o}}^o \|_{h_{o}}^2 
	= -\langle \mathcal{D}^{e\to o} (\hat{A}u^e)_{h_{e}}, \mathcal{P}^o u_{h_{o}}^o \rangle 
	= -\langle \left(\mathcal{P}^o\right)^{-1} \mathcal{Q}^o (\hat{A}u^e)_{h_{e}}, \mathcal{P}^o u_{h_{o}}^o \rangle
	= -\langle \mathcal{Q}^o (\hat{A}u^e)_{h_{e}}, u_{h_{o}}^o \rangle,
\end{align*}
analogously $\frac{1}{2} \partial_t \| u_{h_{e}}^e \|_{h_{e}}^2 = -\langle \mathcal{Q}^e (\hat{A}^T u^o)_{h_{o}}, u_{h_{e}}^e \rangle$ and with the SBP-property
\begin{nalign}
	\label{eqn:discreteEnergyDerivative01}
	\frac{1}{2}\partial_t \left( \| u_h^o \|_{h_{o}}^2 + \| u_h^e \|_{h_{e}}^2 \right) 
	&= -\langle \left(\left.\mathcal{Q}^o\right.^T + \mathcal{Q}^e \right) u_{h_{o}}^o, (\hat{A}u^e)_{h_{e}} \rangle\\
	&= \left(u_{h_{o},0}^o\right)^T\hat{A}u_{h_{e},0}^e-\left(u_{h_{o},N}^o\right)^T\hat{A}u_{h_{e},N}^e.
\end{nalign}
The SAT method consists of augmenting the equations of the boundary nodes by a penalty term that imposes the boundary condition in a weak sense ($a\in\lbrace o,e \rbrace$)
\begin{align}
\partial_t u_{h_{a},0}^a + \ldots &= (P^a_{0,0})^{-1} \tau^a_0 \left(u^o_{h_{o},0} - \left( -L \hat{A} u^e_{h_{e},0} + g_{h,0} \right)\right) 
\label{eqn:SBPFDSAT1Dupdate1}
\\
\partial_t u_{h_{a},N}^a + \ldots &= (P^a_{N,N})^{-1} \tau^a_N \left(u^o_{h_{o},N} - \left(L \hat{A} u^e_{h_{e},N} + g_{h,N}\right)\right) ,
\label{eqn:SBPFDSAT1Dupdate2}
\end{align} 
with penalization matrices $\tau^o_0,\tau^o_N\in\R^{n_o\times n_0}$ and $\tau^e_0,\tau^e_N\in\R^{n_e\times n_o}$. With theses terms \cref{eqn:discreteEnergyDerivative01} becomes
\newcommand{\uh}[1]{\begin{pmatrix}	u^{o}_{h_{o},#1} \\	u^{e}_{h_{e},#1} \end{pmatrix}}
\begin{nalign}
	\label{eqn:discreteEnergyDerivative02}
	\frac{1}{2}\partial_t \left( \| u_{h_{o}}^o \|_{h_{o}}^2 + \| u_{h_{e}}^e \|_{h_{e}}^2 \right) 
	&= 	
	\uh{0}^T
	\begin{pmatrix}
		 \tau_0^{o} & \left(\tau_0^{e}\right)^T + \hat{A} + \tau_0^{o} L\hat{A} \\
		 0& \tau_0^{e} L\hat{A}
	\end{pmatrix}
	\uh{0}
	-
	\uh{0}^T
	\begin{pmatrix}
		\tau_0^{o} \\
		\tau_0^{e}
	\end{pmatrix}
	g_{h,0}\\
	&+
	\uh{N}^T
	\begin{pmatrix}
		\tau_N^{o} & \left(\tau_N^{e}\right)^T -  \hat{A} - \tau_N^{o} L\hat{A}  \\
		0& -\tau_N^{e} L\hat{A}
	\end{pmatrix}
	\uh{N}
	-
	\uh{N}^T
	\begin{pmatrix}
		\tau_N^{o} \\
		\tau_N^{e}
	\end{pmatrix}
	g_{h,N}.
\end{nalign}
In order to obtain an energy bound on the semi-discrete level, we need to choose the penalization matrices $\tau^o_\ast$, $\tau^e_\ast$ appropriately. We first reduce the complexity by relating $\tau^e_\ast$ to $\tau^o_\ast$ such that the bilinear terms in \cref{eqn:discreteEnergyDerivative02} vanish, i.e.
\begin{equation}
	\tau^{e}_0 = -\left(\hat{A} + \tau^{o}_0 L \hat{A} \right)^T
	\quad\text{and}\quad 
	\tau^{e}_N = \left(\hat{A} + \tau^{o}_N L \hat{A} \right)^T.
	\label{eqn:stabilityConditionsSAT01}
\end{equation}
With \cref{eqn:stabilityConditionsSAT01} we obtain necessary and sufficient conditions for $\tau^o_0$, $\tau^o_N$ to assure an energy bound analogous to the continuous version \eqref{eqn:energyStabilityRT} for the semi-discrete approximation.
\begin{theorem}[Discrete Energy Bound (1D)]
	A solution to \cref{prob:1dAdvection} discretized in space using SBP-FD and SAT as described above satisfies the energy bound
	\begin{equation}
		\| u_{h_{o}}^o(T) \|_{h_{o}}^2 + \| u_{h_{e}}^e(T) \|_{h_{e}}^2 
		\leq 
		\| (u^o_0)_{h_{o}} \|_{h_{o}}^2 + \| (u^e_0)_{h_{e}} \|_{h_{e}}^2 
		+ C \left(\| g_L \|^2_{L^2(0,T)} + \| g_R \|^2_{L^2(0,T)} \right)
		\label{eqn:discreteEnergyBound1D}
	\end{equation}
	with $C=\max  \bigcup\limits_{\star\in\lbrace 0,N \rbrace} \lbrace  \| \tau^{o}_\star \|_2, \| L^+ + {\tau_\star^o}^T \|_2 \rbrace $, \\iff the penalization matrices $\tau^o_\star$ ($\star\in\lbrace 0,N \rbrace$) are semi-negative definite and $\forall x\in\R^{n_o}$
	\begin{equation}
		x^T L \left(-\tau_\star^o\right)^T L x \leq x^T L x.
		\label{eqn:discreteEnergyBound1DCondition}
	\end{equation}
\label{thrm:energyBoundSBPSAT1D}	
\end{theorem}
\begin{proof}
	Assume the energy bound \eqref{eqn:discreteEnergyBound1D} holds.
	For the energy bound to hold for all $u^{o}_{h_{o},\star}$ and $u^{e}_{h_{e},\star}$ at all times $T$ for the particular case of vanishing boundary source ($g=0$), we can deduce from \cref{eqn:discreteEnergyDerivative02} that the matrices $\tau_\star^o$, $( \tau_0^{e} L \hat{A} ) =  - ((\hat{A} + \tau^{o}_0 L \hat{A} )^T L \hat{A}$ and $-(\tau^{e}_N L \hat{A}) = - (\hat{A} + \tau^{o}_N L \hat{A} )^T L \hat{A} $ have to be semi-negative definite. Condition \eqref{eqn:discreteEnergyBound1DCondition} now follows from
	\begin{equation}
		(\hat{A} + \tau^{o}_\star L \hat{A} )^T L \hat{A} 
		= \hat{A}^T \left( L +  L \left( \tau_\star^o\right )^T L \right) \hat{A} \overset{L\text{ inv.}}{=} \hat{A}^T L \left( L^{-1} +  \left( \tau_\star^o\right )^T \right) L \hat{A},
	\end{equation}
	$\rank{\hat A}=n_o$ and $L$ semi-positive definite.\\
	Now assume that $\tau^o_\star$ is semi-negative definite and \eqref{eqn:discreteEnergyBound1DCondition} holds. 
	We derive a bound of the right hand side of \cref{eqn:discreteEnergyDerivative02} by bounding terms involving $u^{o}_{h_{o},\star}$ and $u^{e}_{h_{e},\star}$ separately. We will use 
	\begin{align}
		x^T M x - 2 x^T M p
		&= (x-p)^T M (x-p) - p^T M p \leq |p^T M p| 
		\leq {\|M\|}_2 {\|p\|}_2^2 
		= {\|M\|}_2\, p^T p
		\label{eqn:boundOnNegative}
	\end{align}
	for semi-negative definite $M\in\R^{n \times n}$.
	With \cref{eqn:boundOnNegative} can conclude the following bound for the terms in $u^{o}_{h_{o},\star}$
	\begin{equation}
		(u^{o}_{h_{o},\star})^T \tau^{o}_\star u^{o}_{h_{o},\star} - (u^{o}_{h_{o},\star})^T \tau^{o}_\star g_\star
		= (u^{o}_{h_{o},\star})^T \tau^{o}_\star u^{o}_{h_{o},\star} - 2(u^{o}_{h_{o},\star})^T \tau^{o}_\star \frac{g_\star}{2}
		\leq \frac{1}{4} {\| \tau^{o}_\star \|}_2 \, {g_\star}^T g_\star.
		\label{eqn:boundOdd}
	\end{equation}
	As $-( L + L \left(\tau_\star^o\right)^T L )$ is semi-negative definite we can bound the terms in $u^{e}_{h_{e},\star}$ similarly 
	\begin{align}
		& (u^{e}_{h_{e},\star})^T (-({\hat A}^T L{\hat A} + {\hat A}^T L \left(\tau_\star^o\right)^T L {\hat A} )) u^{e}_{h_{e},\star} \pm ( u^{e}_{h_{e},\star})^T ({\hat A}^T + {\hat A}^TL(\tau^o_\star)^T)  g_\star \\
		= &(L \hat A u^{e}_{h_{e},\star})^T (-(L^+ + \left(\tau_\star^o\right)^T  )) ( L \hat A u^{e}_{h_{e},\star}) \pm ( L \hat A u^{e}_{h_{e},\star})^T (L^+ + (\tau^o_\star)^T) g_\star \\
		\leq & \frac{1}{4} {\| L^+ + {\tau_\star^o}^T \|}_2 \, {g_\star}^T g_\star.
		\label{eqn:boundEven}
	\end{align}
	The energy bound now follows from
	\begin{equation}
		\frac{1}{2}\partial_t \left( \| u_{h_{o}}^o \|_{h_{o}}^2 + \| u_{h_{e}}^e \|_{h_{e}}^2 \right) 
		\leq \sum_{\star\in\lbrace 0,N \rbrace} \frac{1}{4} \left( \| \tau^{o}_\star \|_2 + \| L^+ + {\tau_\star^o}^T \|_2 \right) \, {g_\star}^T g_\star 
		\leq  \frac{C}{2} \sum_{\star\in\lbrace 0,N \rbrace} {g_\star}^T g_\star
	\end{equation} 
	with $C=\max  \bigcup\limits_{\star\in\lbrace 0,N \rbrace} \lbrace  {\| \tau^{o}_\star \|}_2, {\| L^+ + {\tau_\star^o}^T \|}_2 \rbrace $.
\end{proof}

The most apparent choice for the penalization matrix $\tau_\star^o$ to assure the discrete energy bound is $\tau_\star^o=0$. For the case that $L$ is invertible, the condition (\ref{eqn:discreteEnergyBound1DCondition}) simplifies to 
$$
0 \leq x^T \left(-\tau_\star^o\right)^T x \leq x^T L^{-1} x
$$
and an obvious family of penalization matrices assuring energy stability would be 
$\tau_\star^o=-\alpha L^{-1}$ with $\alpha\in[0,1]$.

\subsection{Extension to 3D}
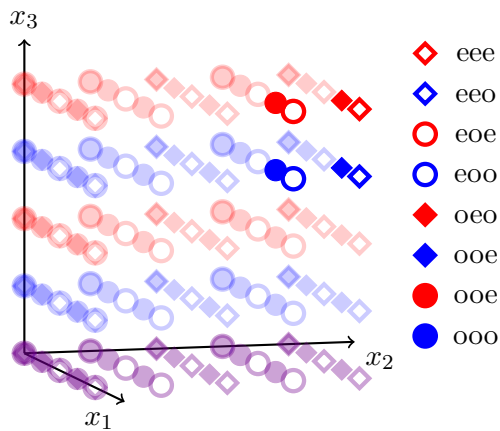
\begin{figure}
	\centering
	\tdplotsetmaincoords{82.5}{75}
\def\symbolsize{8}
\def\innersize{2}
\newcommand{\drawSBPgrid}[6]{
	\foreach \x in #1{
		\foreach \y	in #2{
			\foreach \z	in #3{
				\node[#5,inner sep=\innersize pt,minimum size=\symbolsize pt] at (\x,\y,\z) [draw,#4,#6] {};
				if
			}
		}
	}
}

\begin{tikzpicture}[scale=0.6,tdplot_main_coords]
	
	\coordinate (O) at (0,0,0);
	\def\setEven{0,3,6}
	\def\setOdd{0,1.5,4.5}
	\def\L{7}
	
	\draw[thick,->] (0,0,0) -- (\L+1.5,0,0) node[anchor=north east]{$x_1$};
	\draw[thick,->] (0,0,0) -- (0,\L+0.5,0) node[anchor=north west]{$x_2$};
	\draw[thick,->] (0,0,0) -- (0,0,\L) node[anchor=south]{$x_3$};

	\drawSBPgrid{\setEven}{\setEven}{\setEven}{red}{diamond}{solid,ultra thick,opacity=0.2};
	\drawSBPgrid{\setEven}{\setEven}{\setOdd}{blue}{diamond}{solid,ultra thick,opacity=0.2};
	\drawSBPgrid{\setEven}{\setOdd}{\setEven}{red}{circle}{solid,ultra thick,opacity=0.2};
	\drawSBPgrid{\setEven}{\setOdd}{\setOdd}{blue}{circle}{solid,ultra thick,opacity=0.2};
	
	\drawSBPgrid{\setOdd}{\setEven}{\setEven}{red}{diamond}{fill,opacity=0.2};
	\drawSBPgrid{\setOdd}{\setEven}{\setOdd}{blue}{diamond}{fill,opacity=0.2};
	\drawSBPgrid{\setOdd}{\setOdd}{\setEven}{red}{circle}{fill,opacity=0.2};
	\drawSBPgrid{\setOdd}{\setOdd}{\setOdd}{blue}{circle}{fill,opacity=0.2};
	
	\drawSBPgrid{6}{6}{6}{red}{diamond}{solid,ultra thick,opacity=1.0};
	\drawSBPgrid{6}{6}{{4.5}}{blue}{diamond}{solid,ultra thick,opacity=1.0};
	\drawSBPgrid{6}{{4.5}}{6}{red}{circle}{solid,ultra thick,opacity=1.0};
	\drawSBPgrid{6}{{4.5}}{{4.5}}{blue}{circle}{solid,ultra thick,opacity=1.0};
	
	\drawSBPgrid{{4.5}}{6}{6}{red}{diamond}{fill,opacity=1.0};
	\drawSBPgrid{{4.5}}{6}{{4.5}}{blue}{diamond}{fill,opacity=1.0};
	\drawSBPgrid{{4.5}}{{4.5}}{6}{red}{circle}{fill,opacity=1.0};
	\drawSBPgrid{{4.5}}{{4.5}}{{4.5}}{blue}{circle}{fill,opacity=1.0};
	
	\def\xlegend{0}
	\def\ylegend{1.3*\L}
	\def\zlegend{0.9*\L}
	\def\deltaZ{0.9}
	\def\deltaY{1.15}
	\node[diamond,inner sep=\innersize pt,minimum size=\symbolsize pt] at (\xlegend,\ylegend,\zlegend+0.075) [draw,red,ultra thick,solid] {};
	\node at (\xlegend,\ylegend+\deltaY,\zlegend) {eee};
	\node[diamond,inner sep=\innersize pt,minimum size=\symbolsize pt] at (\xlegend,\ylegend,\zlegend+0.075-1*\deltaZ) [draw,blue,ultra thick] {};
	\node at (\xlegend,\ylegend+\deltaY,\zlegend-1*\deltaZ) {eeo};
	\node[circle,inner sep=\innersize pt,minimum size=\symbolsize pt] at (\xlegend,\ylegend,\zlegend+0.075-2*\deltaZ) [draw,red,ultra thick] {};
	\node at (\xlegend,\ylegend+\deltaY,\zlegend-2*\deltaZ) {eoe};
	\node[circle,inner sep=\innersize pt,minimum size=\symbolsize pt] at (\xlegend,\ylegend,\zlegend+0.075-3*\deltaZ) [draw,blue,ultra thick] {};
	\node at (\xlegend,\ylegend+\deltaY,\zlegend-3*\deltaZ) {eoo};
	\node[diamond,inner sep=\innersize pt,minimum size=\symbolsize pt] at (\xlegend,\ylegend,\zlegend+0.075-4*\deltaZ) [draw,red,ultra thick,fill] {};
	\node at (\xlegend,\ylegend+\deltaY,\zlegend-4*\deltaZ) {oeo};
	\node[diamond,inner sep=\innersize pt,minimum size=\symbolsize pt] at (\xlegend,\ylegend,\zlegend+0.075-5*\deltaZ) [draw,blue,ultra thick,fill] {};
	\node at (\xlegend,\ylegend+\deltaY,\zlegend-5*\deltaZ) {ooe};
	\node[circle,inner sep=\innersize pt,minimum size=\symbolsize pt] at (\xlegend,\ylegend,\zlegend+0.075-6*\deltaZ) [draw,red,ultra thick,fill] {};
	\node at (\xlegend,\ylegend+\deltaY,\zlegend-6*\deltaZ) {ooe};
	\node[circle,inner sep=\innersize pt,minimum size=\symbolsize pt] at (\xlegend,\ylegend,\zlegend+0.075-7*\deltaZ) [draw,blue,ultra thick,fill] {};
	\node at (\xlegend,\ylegend+\deltaY,\zlegend-7*\deltaZ) {ooo};
\end{tikzpicture}
	\caption{Illustration of staggered grids in three-dimensions. The grid boundary is assumed to be on the surfaces spanned by the axes. Each of the eight symbols stores a certain odd/even class of variables. In total there are eight classes. In the interior there is the following pattern: In $x_1$-direction filled and open variable classes interchange, in $x_2$-direction diamonds and circles interchange, in $x_3$-direction blue and red variable classes interchange. Analogously to the 1d setting in Fig.~\ref{fig:staggeredGrid1D} the respective variable classes are jointly present on the boundary surfaces.}
	\label{fig:staggeredGrid}
\end{figure}

Before explaining the extension to the three-dimensional case, we realize that the even/odd classification with respect to the three cartesian coordinate axes divides the spherical harmonic functions into $8\ (=2^3)$ families with two families coupling within the recursion relation of direction $i$, iff their type is distinct only in direction $i$. We classify the respective expansion coefficients accordingly. This insight into the coupling allows us to apply the staggered SBP-FD discretization to the three-dimensional \PN{N} equations.
We discretize the eight families of expansion coefficients on an rectangular domain $G=[x_{1,L},x_{1,R}]\times[x_{2,L},x_{2,R}]\times[x_{3,L},x_{3,R}]\subset\R^3$ using 8 staggered grids, such that the coupling structure within the equations is respected, see Figure \ref{fig:staggeredGrid}. 
When representing grid functions by vectors of node values with $x_3$ as the contiguous direction, following $x_2$ and then the $x_1$ direction, the matrices identifying the staggered SBP difference operators in three dimensions can be written using the Kronecker product
\def\grididx{a}
\def\gridid{\bold{a}}
\begin{equation}
	D^{\gridid}_1 = D^{\grididx_1}_{x_1} \otimes I_{x_2}^{\grididx_2} \otimes I_{x_3}^{\grididx_3},
	\quad
	D^{\gridid}_2 = I_{x_1}^{\grididx_1} \otimes D^{\grididx_2}_{x_2} \otimes I_{x_3}^{\grididx_3},
	\quad
	D^{\gridid}_3 = I_{x_1}^{\grididx_1} \otimes I_{x_2}^{\grididx_3} \otimes D^{\grididx_3}_{x_3}.
	\label{eqn:3DSBPFD}
\end{equation}
The triplet $\gridid=(a_1,a_2,a_3)$ with $\grididx_i\in\lbrace o,e \rbrace$ identifies grids $h_\gridid$ and variables $u_{h_{\gridid}}^\gridid$ with respect to the odd-even classification in three dimensional space and $I_{x_1}^o$, etc. denote identity matrices of the same dimension as the number of grid nodes in the direction indicated by the subscript. The same approach can be used to construct scalar products assigned to the different grids
\begin{align}
	{\langle f_{h_{\gridid}}, g_{h_{\gridid}} \rangle}_{h_\gridid} \eqdef \langle f_{h_{\gridid}}, 
		(\mathcal{P}_{x_1}^{\grididx_1} \otimes \mathcal{P}_{x_{2}}^{\grididx_2} \otimes \mathcal{P}_{x_{3}}^{\grididx_3})
	\ g_{h_{\gridid}} \rangle 
	\qquad
	\| f_{h_\gridid} \|^2_{h_\gridid} = {\langle f_{h_\gridid}, f_{h_\gridid} \rangle}_{h_\gridid}
	.
\end{align}
and for the boundary grids $\partial h_{\gridid}=h_{\gridid}\cap \overline{G}$
\begin{align}
	{\langle f_{\partial h_{\gridid}}, g_{\partial h_{\gridid}} \rangle}_{\partial h_{\gridid}}
	= \sum_{i\in\partial h_{\gridid}} w_i^\gridid
		\left(f_{\partial h_{\gridid}}\right)^T_i \left( g_{\partial h_{\gridid}}\right)_i
	\qquad
	\| f_{\partial h_{\gridid}} \|^2_{\partial h_{\gridid}} = {\langle f_{\partial h_{\gridid}}, g_{\partial h_{\gridid}} \rangle}_{\partial h_{\gridid}}
\end{align}
with $w^\grididx_i=\sum_{i=1}^d \sum_{\star\in\lbrace 0,N_{d} \rbrace} \delta_{i_d,\star} \prod_{j\neq d} P_{x_j}^{\grididx_j}$. Analogous to the one-dimensional case we discretize in space by replacing the spatial derivatives by the difference operators \eqref{eqn:3DSBPFD} and augmenting the equations of the boundary nodes with SAT terms. Using $c_d(\gridid)$ to denote the complement of odd-even class $\gridid$ in direction $d$ (e.g. $c_2(ooo)=oeo$), and $o_d(\gridid)$ and ($e_d(\gridid)$) the triplet obtained by replacing the $d$-th component of $\gridid$ by $o$ and $e$ respectively (e.g $e_3(\gridid)=(\grididx_1\grididx_2 e)$) we can write the semi-discrete system by
\begin{align}
	\partial_t u_{h_{\gridid},i}^\gridid
	+ \sum_{d=1}^3 
	A^{\gridid}_d 
	\left( \mathcal{D}^{\gridid}_d u^{c_{d}(\gridid)}_{h_{c_{d}(\gridid)}} \right)_i 
	 &= 
	\sum_{\substack{d=1,\ldots,3\\ \star\in\lbrace 0,N_d \rbrace}}
	\frac{\delta_{i_d,\star}}{{(P_{x_{d}}^{\grididx_d})}_{\star,\star}} \tau^{d,\gridid}_\star \left(u^{o_d(\gridid)}_{h_{o_d(\gridid)},i} - \left(\pm L^{e_d(\gridid)}_{d} \hat{A}_d^{e_d(\gridid)} u^{e_d(\gridid)}_{h_{e_d(\gridid)},i} + {\left( g^\grididx_d \right)}_{i} \right)\right),
	\label{eqn:SBPFDSATfull}
\end{align}
where $i$ is a triplet identifying a node of grid $\gridid$. We used $\pm$ to merge the formulas for $i_d=0$ with a minus sign and $i_d=N_d$ with a plus sign. Although blurred by the indices and superscripts the discretization \eqref{eqn:SBPFDSATfull} is the three dimensional counterpart of \cref{eqn:SBPFDSAT1Dupdate1,eqn:SBPFDSAT1Dupdate2} and we can use the same arguments as used in the proof of \cref{thrm:energyBoundSBPSAT1D} to obtain an analogous energy bound, given that 
\begin{equation}
	\tau^{d,e_d(\gridid)}_0 = -\left(\hat{A}_d^{e_d(\gridid)} + \tau^{d,o_d(\gridid)}_0 L_d^{e_d(\gridid)} \hat{A}_d^{e_d(\gridid)} \right)^T
	\quad\text{and}\quad 
	\tau^{d,e_d(\gridid)}_{N_d} = \left(\hat{A}_d^{e_d(\gridid)} + \tau^{d,o_d(\gridid)}_{N_d} L_d^{e_d(\gridid)} \hat{A}_d^{e_d(\gridid)} \right)^T .
\end{equation}
\begin{theorem}[Discrete Energy Bound (3D)]
	A solution to \cref{prob:1dAdvection} discretized in space using SBP-FD and SAT as described above satisfies the energy bound
	\begin{equation}
		\sum_{\gridid\in\lbrace e,o\rbrace^3} {\| u^\grididx_{h_\gridid}(T) \|}_{h_\gridid}^2 \leq 
		\sum_{\gridid\in\lbrace e,o\rbrace^3} {\| {(u^\grididx_0)}_{h_\gridid} \|}_{h_\gridid}^2
		+ C {\| g^\grididx_d \|}^2_{\partial h_\gridid (0,T)}
		\label{eqn:discreteEnergyBound3D}
	\end{equation}
	with $C=\max  \bigcup\limits_{\substack{d=1,\ldots,3\\ \gridid\in\lbrace e,o\rbrace^3}} \lbrace  {\| \tau^{d,o_d(\gridid)} \|}_2, {\| {L_d}^+ + \left( \tau^{d,o_d(\gridid)}\right)^T \|}_2 \rbrace $ and ${\|\, \cdot\, \|}^2_{\partial h_\gridid (0,T)}=\int_0^T \|\, \cdot\, \|_{\partial h_\gridid}^2 dt$, \\
	iff the penalization matrices $\tau^{d,\gridid}$ are semi-negative definite and $\forall x\in\R^{n_\gridid}$
	\begin{equation}
		x^T L_d^{e_d(\gridid)} \left(-\tau^{d,\gridid}\right)^T L_d^{e_d(\gridid)} x \leq x^T L_d^{e_d(\gridid)} x.
		\label{eqn:discreteEnergyBound3DCondition}
	\end{equation}
	\label{thrm:energyBoundSBPSAT3D}	
\end{theorem}

 	\section{Numerical Results}
\label{sec:numericalResults}
\newcommand{\figsizeSub}{0.9\hsize}
\newcommand{\figsizeSingle}{0.45\hsize}
\newcommand{\fontsizeFig}{\LARGE}
\newcommand{\fontsizeSub}{\LARGE}
In this section we present numerical results of four test cases with different objectives. The first test case confirms and visualizes the decay of energy over time for vacuum boundary conditions. The second demonstrates the effect of incoming waves depending on standing waves and the interpretation as artificial boundary sources. Test cases three and four are motivated by electron probe microscopy (EPMA)\cite{Reimer:1998} and aim to address the accuracy of the boundary conditions in an application of engineering relevance. We consider the two boundary scenarios relevant in EPMA: Electrons leaving a solid material into vacuum and beam electrons entering a solid material. For both scenarios we compare different approximation orders $N$ and address the accuracy of $P_N$ solutions through comparison to Monte Carlo (MC) results obtained by taking averages over a set of $10^7$ electron trajectories simulated with the software package DTSA-II\cite{Ritchie:2005,Ritchie:2009}. It is generally accepted that the odd order \PN{N} approximation is superior to the next even order approximation \PN{N+1} and hence, excluding the second test case, we follow the common practice of choosing $N$ odd.\cite{HOWLET:1959}\\
Unless stated otherwise, the presented results were obtained from discretization of \PN{N} equations as explained in \cref{sec:numericalSolution} with second order SBP operators based on \texttt{github.com/ooreilly/sbp} \cite{Oreilly:2017} and with a second order Strang splitting for numerical integration in time.

\subsection{Test Case 1: Energy Stability}

\begin{figure*}[h]
	\centering
	\begin{subfigure}{0.49\hsize}
		\centering
		\includegraphics[width=0.9\hsize]{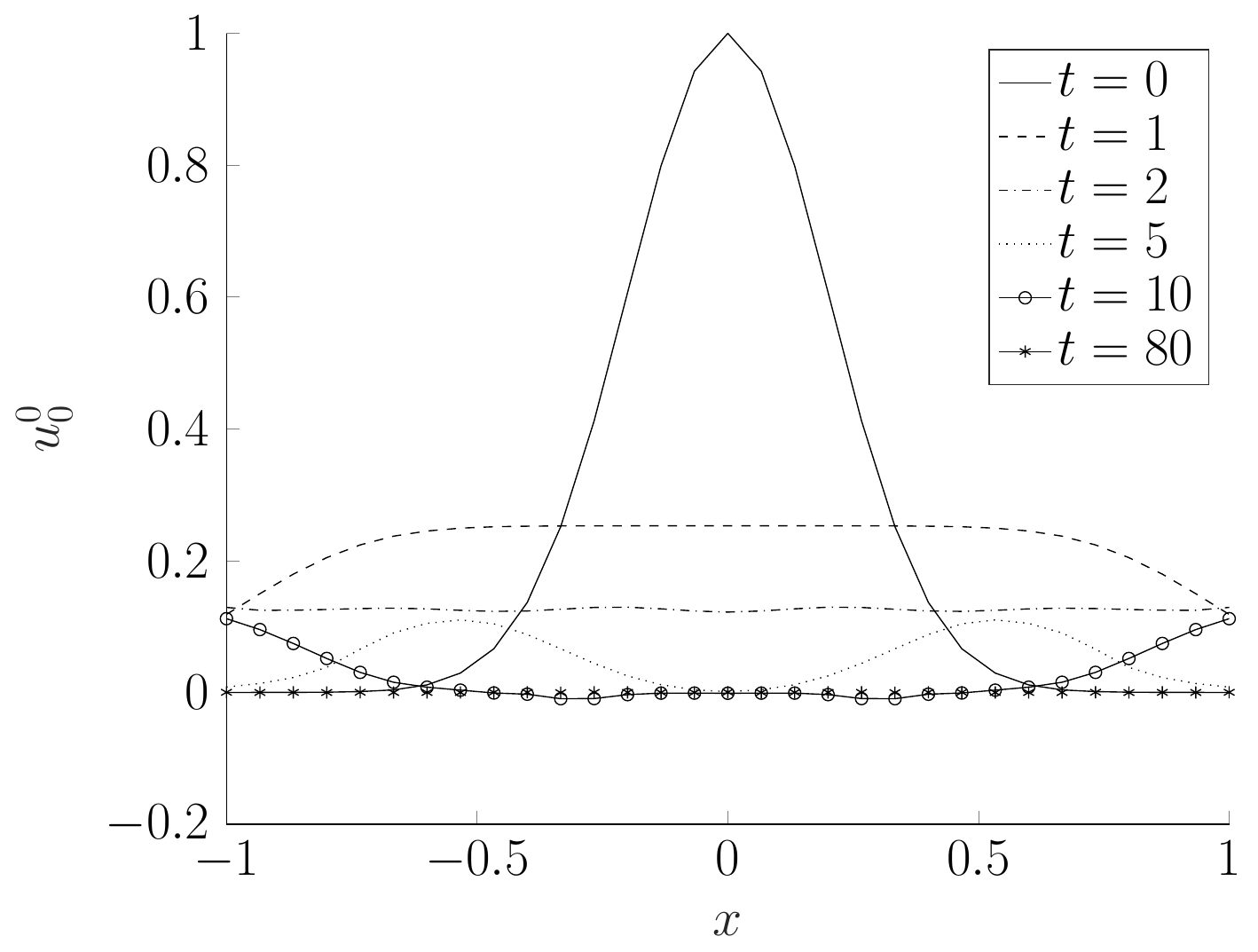}
		\subcaption{Snapshots of a particle density $u_0^0$.}
		\label{fig:tc1evolution}
	\end{subfigure}
	\begin{subfigure}{0.49\hsize}
		\centering
		\includegraphics[width=0.9\hsize]{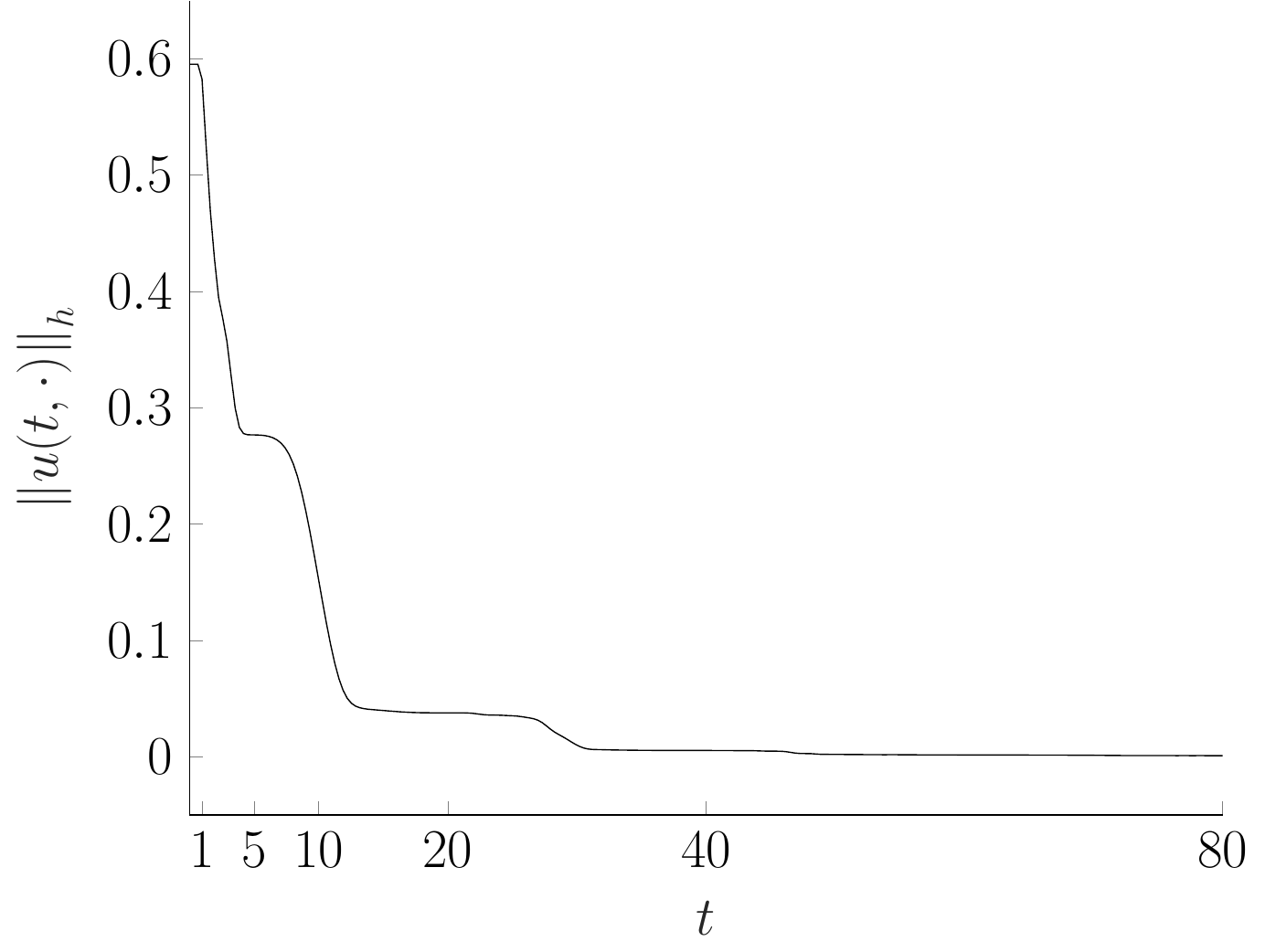}
		\subcaption{Decrease of energy over time.}
		\label{fig:stability}
	\end{subfigure}
	\caption{Results obtained with \PN{13} for an initial bulk of particles with constant distribution in direction space in the absence of collisions and physical boundaries.}
	\label{fig:kernelIssue}
\end{figure*}

In the first test case we consider particles moving in vacuum, i.e. particles are in free-flight and do not interact with a background medium. 
The initial particle distribution is constant in the directions $x_2$ and $x_3$, follows a normal distribution in $x_1$ direction centered at $x_1=0$ and is independent of the direction $\Omega$. We restrict our view onto the spatial domain $(-1,1)\times\R\times\R$ and assume, that the initial particle density is negligible outside, so particles do not enter the domain ($\psi_{in}=0$). 
The $P_N$ model for this case reduces to a pseudo one dimensional problem and, identifying $x_1$ with $x$, we can write the initial condition for the \PN{N} approximation by 
\begin{equation*}
	u_l^k(0,x) = \begin{cases}
		\normpdf{x} & l=k=0 \\
		0 & \text{else}.
	\end{cases}
\end{equation*}
\Cref{fig:tc1evolution} shows the particle density $u_0^0$ at different times $t^*$ for the initial distribution with $\mu_x=0$ and $\sigma_x = 0.2$ obtained with approximation order $N=13$. The density profile smooths out symmetrically towards both boundaries and approaches zero, which coincides with the picture of the underlying particle system, where particles leave the domain and no new particles enter over the boundary. The solution also illustrates two artifacts of the \PN{N} method: The density is not assured to be positive and the set of velocities at which information is transported is finite.\\
In \Cref{fig:stability} the discrete SBP energy norm of the solution is plotted over time. The graph confirms the decrease of the energy in time for vacuum boundary conditions. As the set of transport velocities is finite, the initial data centered around $x=0$ reach the boundaries in a package manner, which leads to the terraced shape of the graph.

\subsection{Test Case 2: Artificial Boundary Source}

\begin{figure}[h]
	\centering
	\includegraphics[width=0.45\hsize]{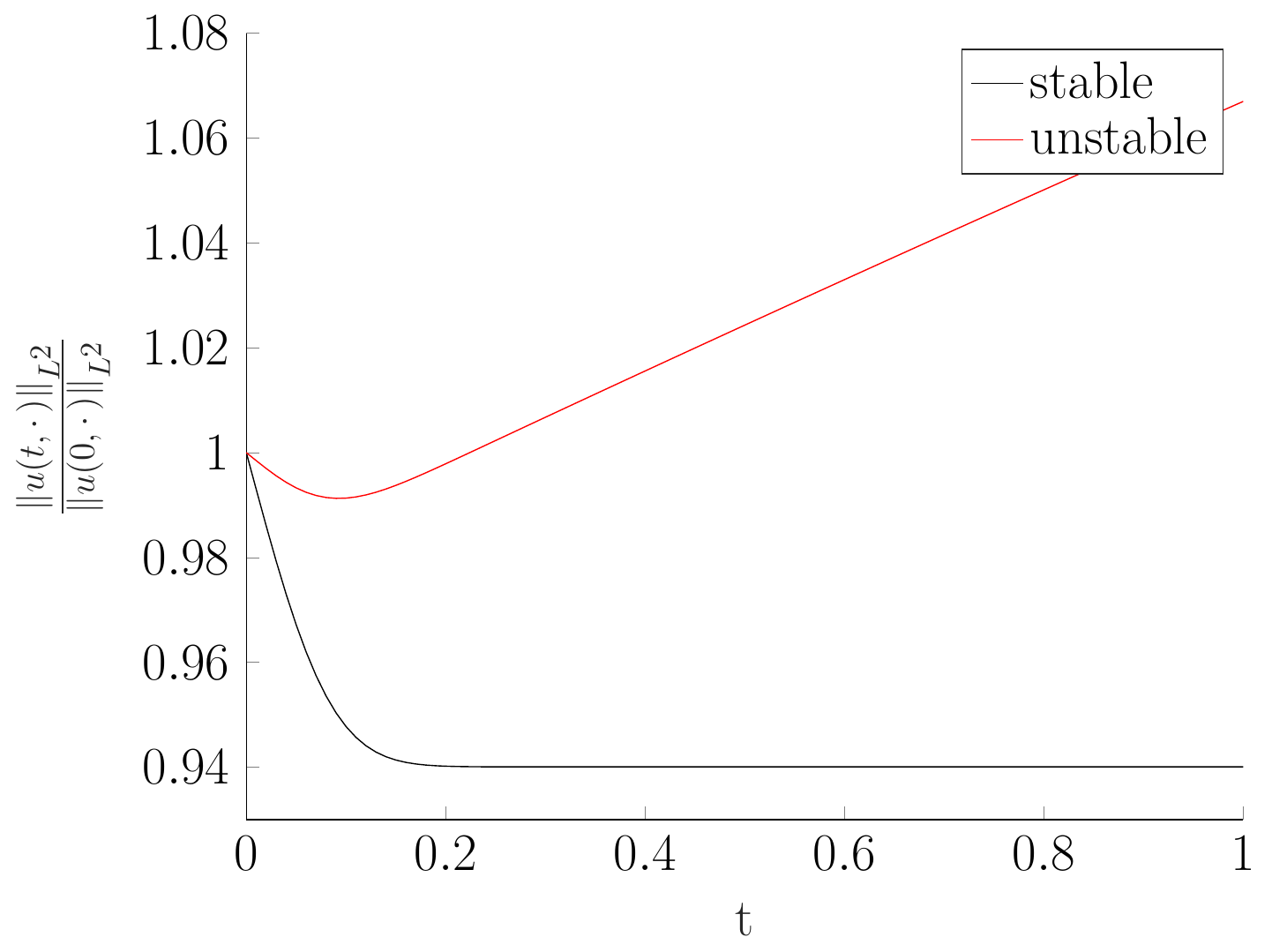}
	\caption{Evolution of energy for test case 2 with stable and unstable boundary conditions; Increase of energy in the absence of boundary source terms caused by standing waves acting as artificial source terms in unstable boundary condition.}
	\label{fig:kernelIssueNorm}
\end{figure}

The purpose of this case is to demonstrate the necessity of stable boundary conditions, to assure a strict energy bound, i.e. the energy does not increase in the absence of boundary sources ($g=0$). As for the first test case, we neglect scattering, consider pseudo one dimensional \PN{N} equations ($N=2$) on the $x_1$ axis and identify $x_1$, $x_2$ and $x_3$ by $x$, $y$ and $z$. We choose the domain to be $G=\R_{+}$, and construct an initial condition that reduces the number of non-trivial coefficients to only four and allows to consider solely the domain close to the boundary at $x=0$. We set all coefficients $u_l^k$ to zero, except for the four coefficients that are even in both $y$ and $z$ direction, namely $u^{e}=(u_0^0,u_2^0,u_2^2)^T$ and $u^{o}=u_1^1$. The pseudo one dimensional \PN{2} equations are
\begin{equation*}
	\partial_t \begin{pmatrix}
	u^o \\ u^e
	\end{pmatrix}
	+
	\begin{pmatrix}
	0 & \hat{A} \\ \hat{A}^T & 0
	\end{pmatrix}
	\partial_x
	\begin{pmatrix}
	u^o \\ u^e
	\end{pmatrix}
	= 0, \qquad\text{with }\quad 
	\hat{A}=\begin{pmatrix}	0.5773\ldots & -0.2581\ldots & 0.4472\ldots \end{pmatrix}.
\end{equation*}
For the demonstration we consider the boundary condition with zero source
\begin{equation}
	u^o = \begin{pmatrix} 0.8660\ldots & -0.2420\ldots & 0.4192\ldots \end{pmatrix}
	u^e {\color{gray}\, +\, \cancelto{0}{g}}\ ,
	\label{eqn:demonstrationBC}
\end{equation}
obtained from the submatrix of $\tilde M$ in \cref{eqn:BCwithMcomplete} that couples $u^{o}$ to $u^{e}$. This boundary condition prescribes the incoming waves dependent of standing waves such that no strict stability is given, which can be avoided by neglecting the higher order coefficients of the recursion relation appearing in matrix $\tilde M$ (\cref{eqn:Mcomplete}), as described in \cref{ssec:PNBConsager}. The initial condition is constructed to serve the purpose: We choose $u^{e}$ and compute $u^{o}$ from the boundary condition
\begin{equation*}
	u^{e}|_{t=0} \overset{!}{=} \begin{pmatrix} 1 \\ 2.5 \\ -1 \end{pmatrix} e^{\left(\frac{x}{0.1}\right)^2}
	\quad\overset{BC}{\longrightarrow}\quad
	u^{o}|_{t=0} 
	= -0.1583\ldots\,\, e^{\left(\frac{x}{0.1}\right)^2}.
\end{equation*}
The red graph in \Cref{fig:kernelIssueNorm} shows the change of the energy over time for the unstable boundary condition \eqref{eqn:demonstrationBC}. The change of the energy results from two opposing processes: The energy-loss caused by the outgoing waves, and the energy-gain caused by incoming waves prescribed by the boundary condition at $x=0$. While the energy loss is initially dominant, both processes balance out around $t=0.1$, and as time increases the outgoing waves vanish and we can observe the constant increase in energy caused by the standing waves acting as artificial sources in the boundary condition.\\
The black graph in \Cref{fig:kernelIssueNorm} visualizes the energy evolution for the stabilized version of the test case, i.e. replacing boundary condition \eqref{eqn:demonstrationBC} by the stabilized version as described in \cref{ssec:PNBConsager} and adapting the initial condition of $u^o$. For the stable boundary conditions the energy decays and, as incoming waves are independent of standing waves, approaches steady state when outgoing waves become negligible.


\subsection{Boundary Conditions for EPMA}
We now consider two practically relevant test cases motivated by EPMA, a method to quantify the micro structure of materials\cite{Reimer:1998}. The method is based on exciting a small polished material sample by a focused electron beam within a vacuum chamber. In the context of EPMA the electron transport is usually described in terms of the electron fluence modeled by BCSD \cref{eqn:BCSD}. The electron fluence density, given by the zeroth order moment $u_0^0$ in the series expansion, is the quantity of interest for EPMA and will be used for comparison and evaluation. For the sake of demonstration we consider only two dimensions and assume solutions to be constant in the third direction.\\
In both test cases we consider a sample of pure copper of density 8960 kg/m$^3$. Our \PN{N} code and the MC code DTSA-II implement different models of stopping power and scattering cross-sections. For comparability of the MC and \PN{N} results we set the stopping power to a constant $S=-2\cdot10^{-13}$ (J/m)/(kg/m$^3$) in both codes. The scattering cross-sections used in the \PN{N} simulations were taken from the database of the ICRU report \cite{ICRU77:2007} which was generated using the ELSEPA code\cite{Salvat:2005} and the MC simulations employed scattering cross sections from a NIST database\cite{Jablonski:2003} that was generated by an earlier version of the ELSEPA code.

\subsubsection{Test Case 3: Vacuum}
\begin{figure}
	\centering
	\begin{subfigure}{0.49\hsize}
		\centering
		\includegraphics[width=0.9\hsize]{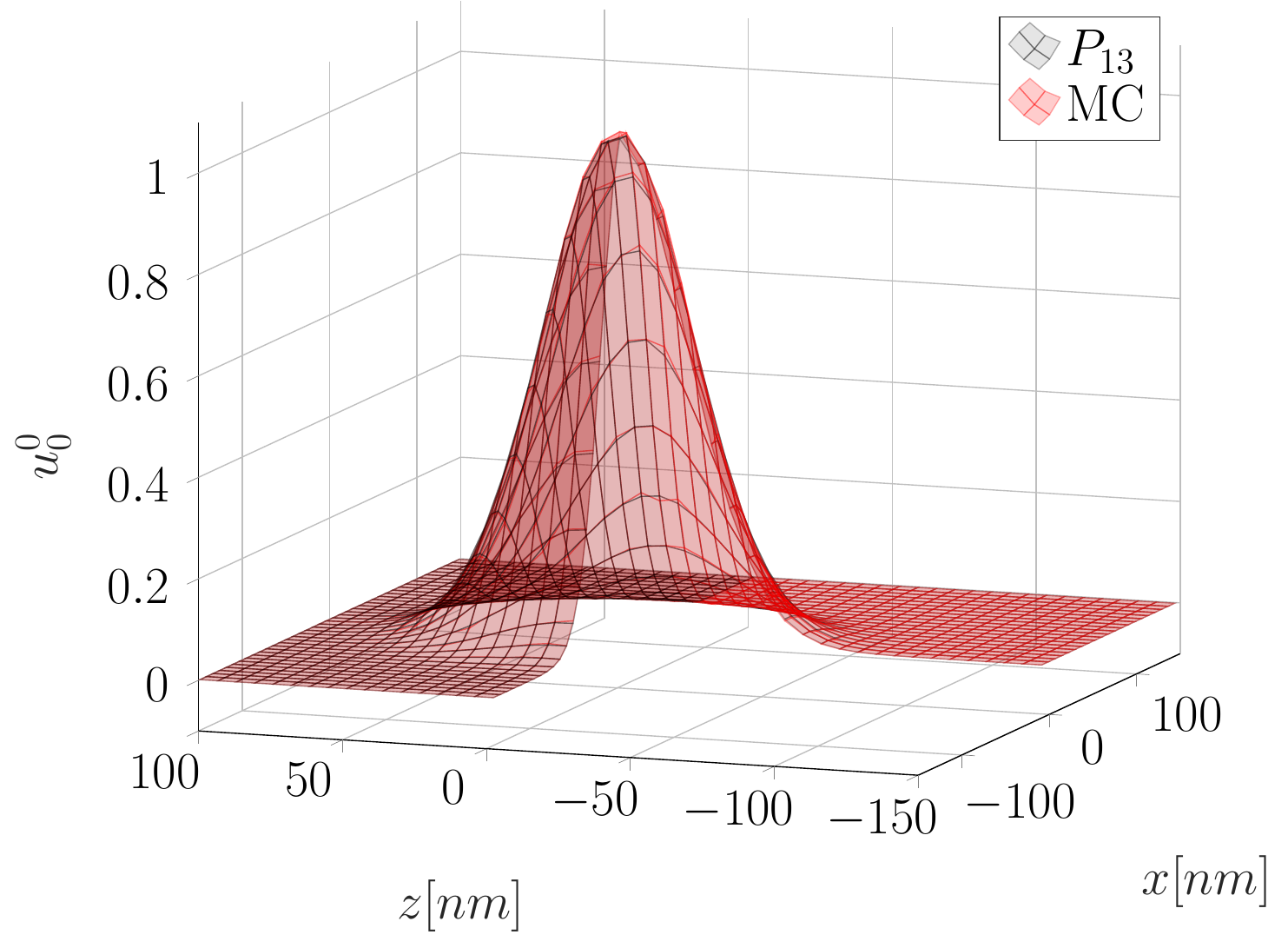}
		\subcaption{$u_0^0$ at $15.0$keV}
		\label{fig:EPMAvacuum2DInitial}
	\end{subfigure}
	\begin{subfigure}{0.49\hsize}
		\centering
		\includegraphics[width=0.9\hsize]{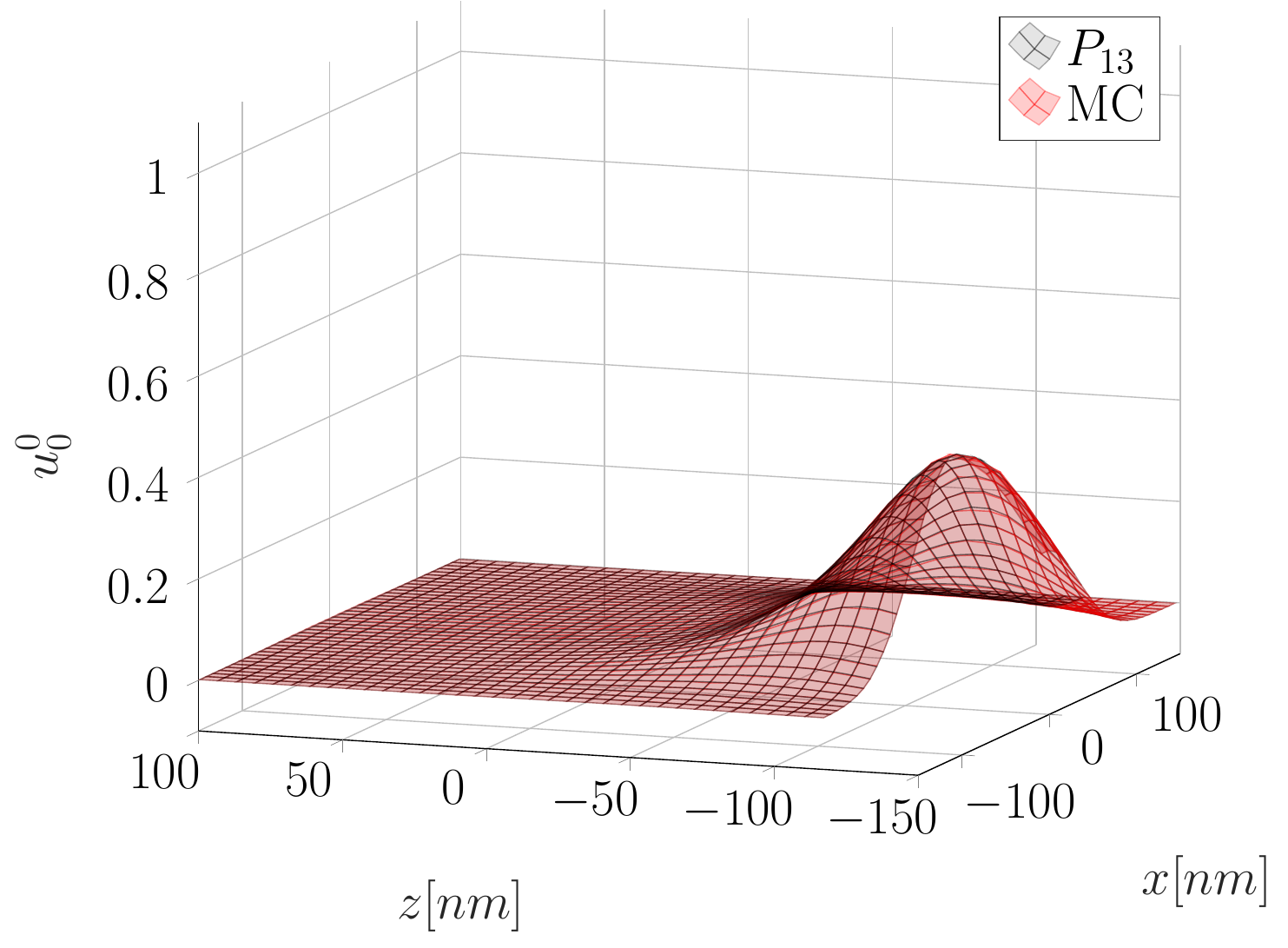}
		\subcaption{$u_0^0$ at $13.5$keV}
	\end{subfigure}
	\begin{subfigure}{0.49\hsize}
		\centering
		\includegraphics[width=0.9\hsize]{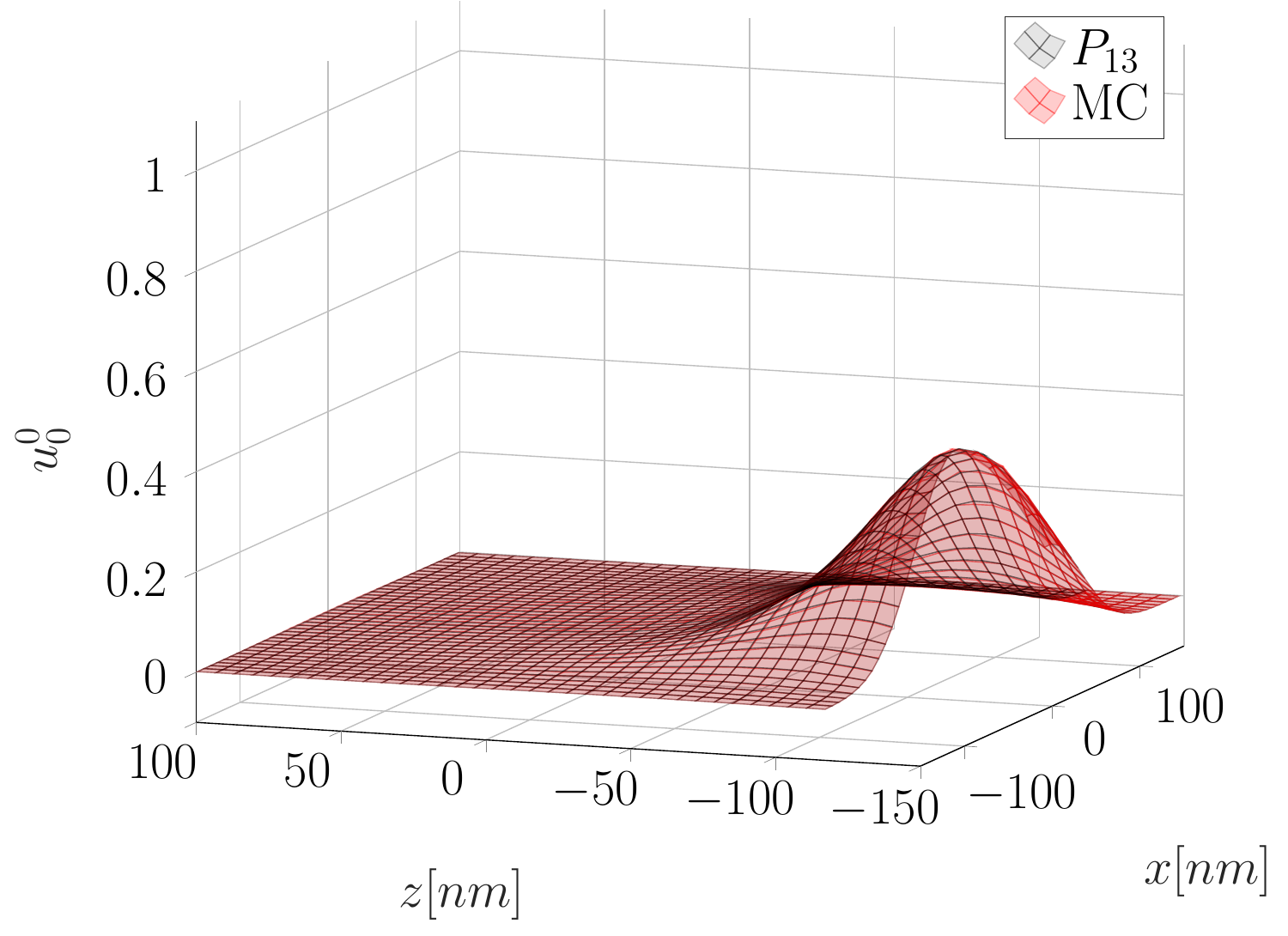}
		\subcaption{$u_0^0$ at $13.0$keV}
	\end{subfigure}
	\begin{subfigure}{0.49\hsize}
		\centering
		\includegraphics[width=0.9\hsize]{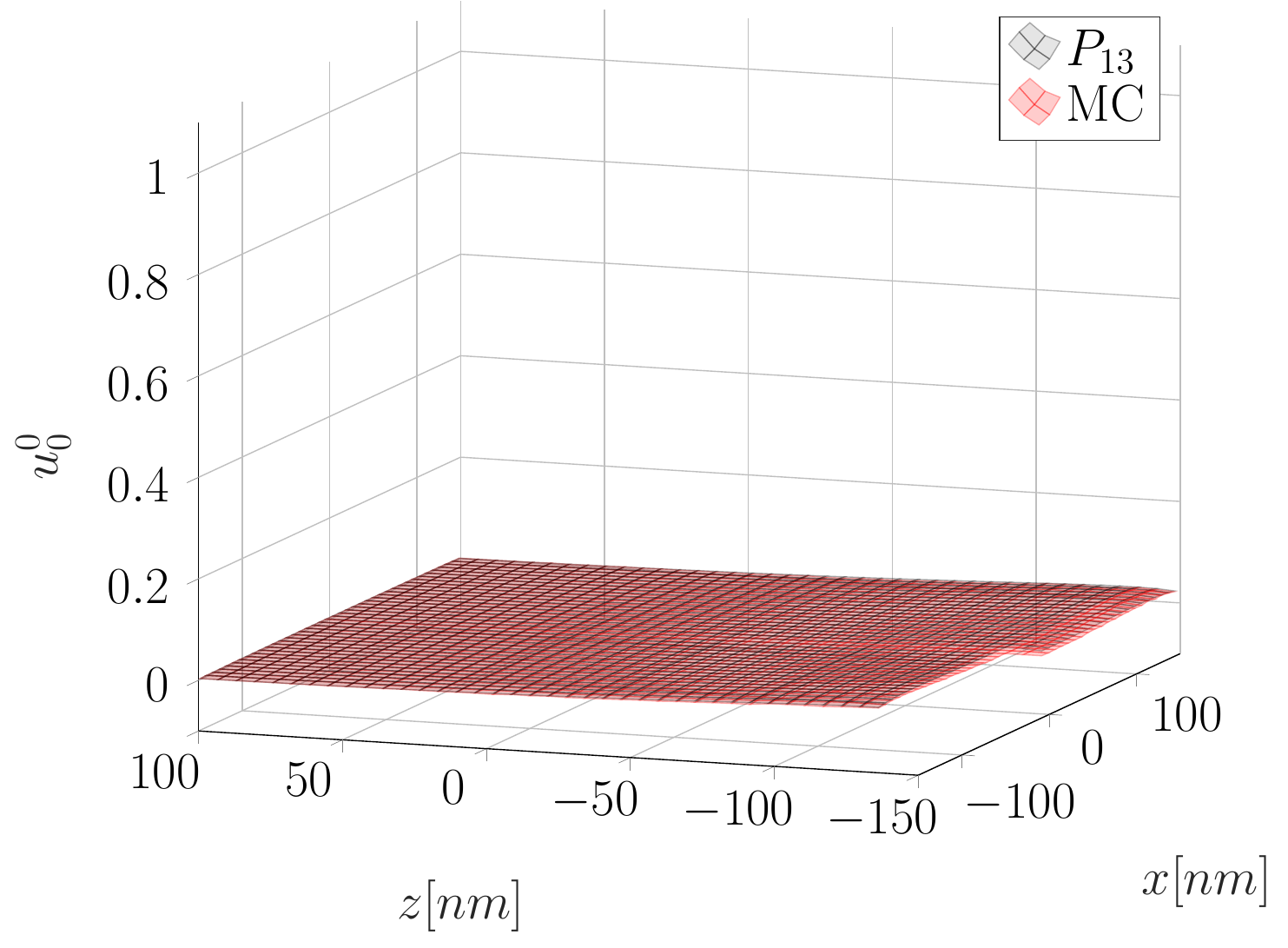}
		\subcaption{$u_0^0$ at $12.0$keV}
	\end{subfigure}
	\caption{Fluence density $u_0^0$ of a bulk of electrons escaping a material sample of pure copper into vacuum over a polished surface at z=$-$150nm computed with \PN{13} and MC.}	
	\label{fig:EPMAvacuum2D}
\end{figure}

\begin{figure}
	\centering
	\begin{subfigure}{0.49\hsize}
		\centering
		\includegraphics[width=0.9\hsize]{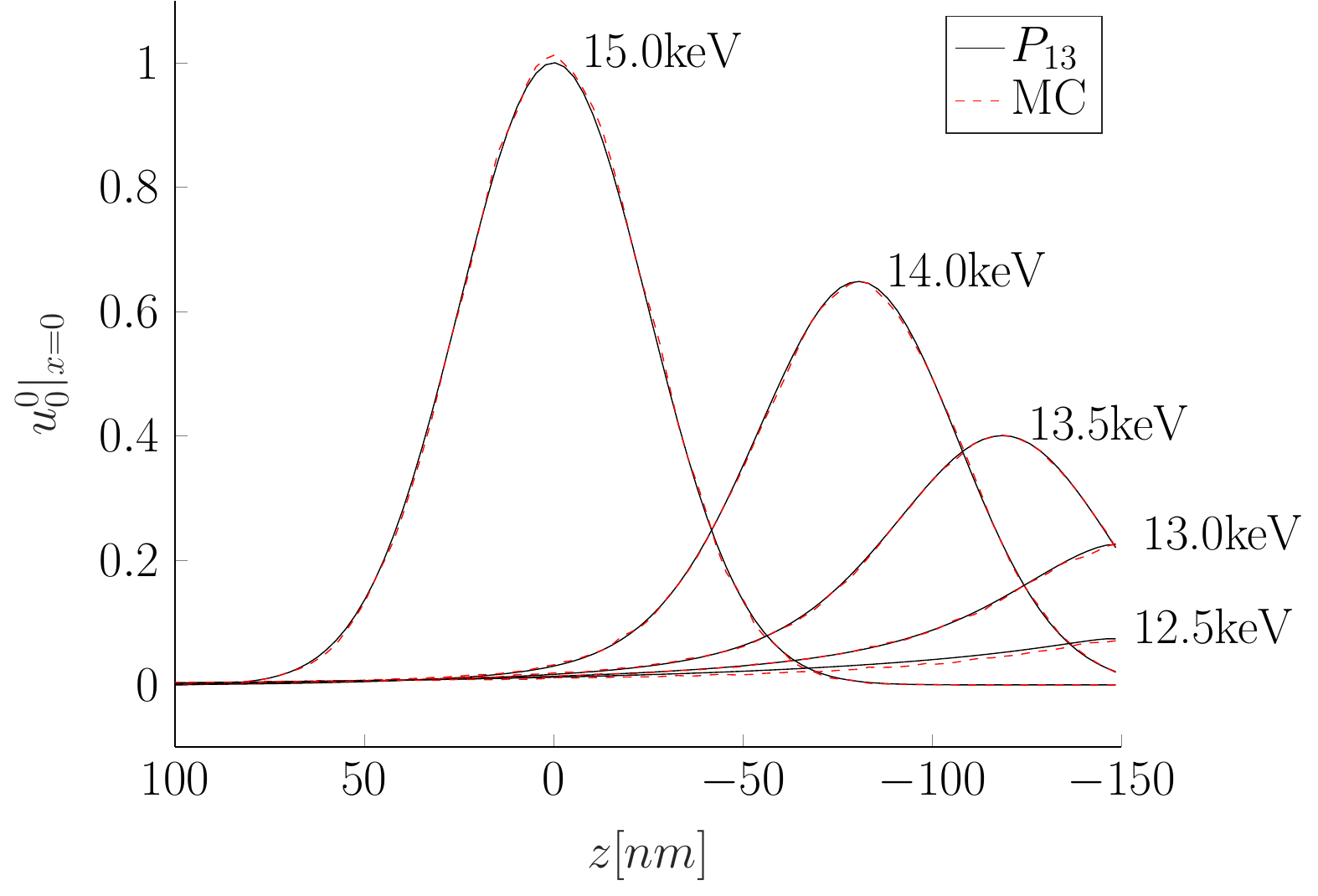}
		\subcaption{\PN{13} vs. MC}
		\label{fig:EPMAvacuumOneD1}
	\end{subfigure}
	\begin{subfigure}{0.49\hsize}
		\centering
		\includegraphics[width=0.9\hsize]{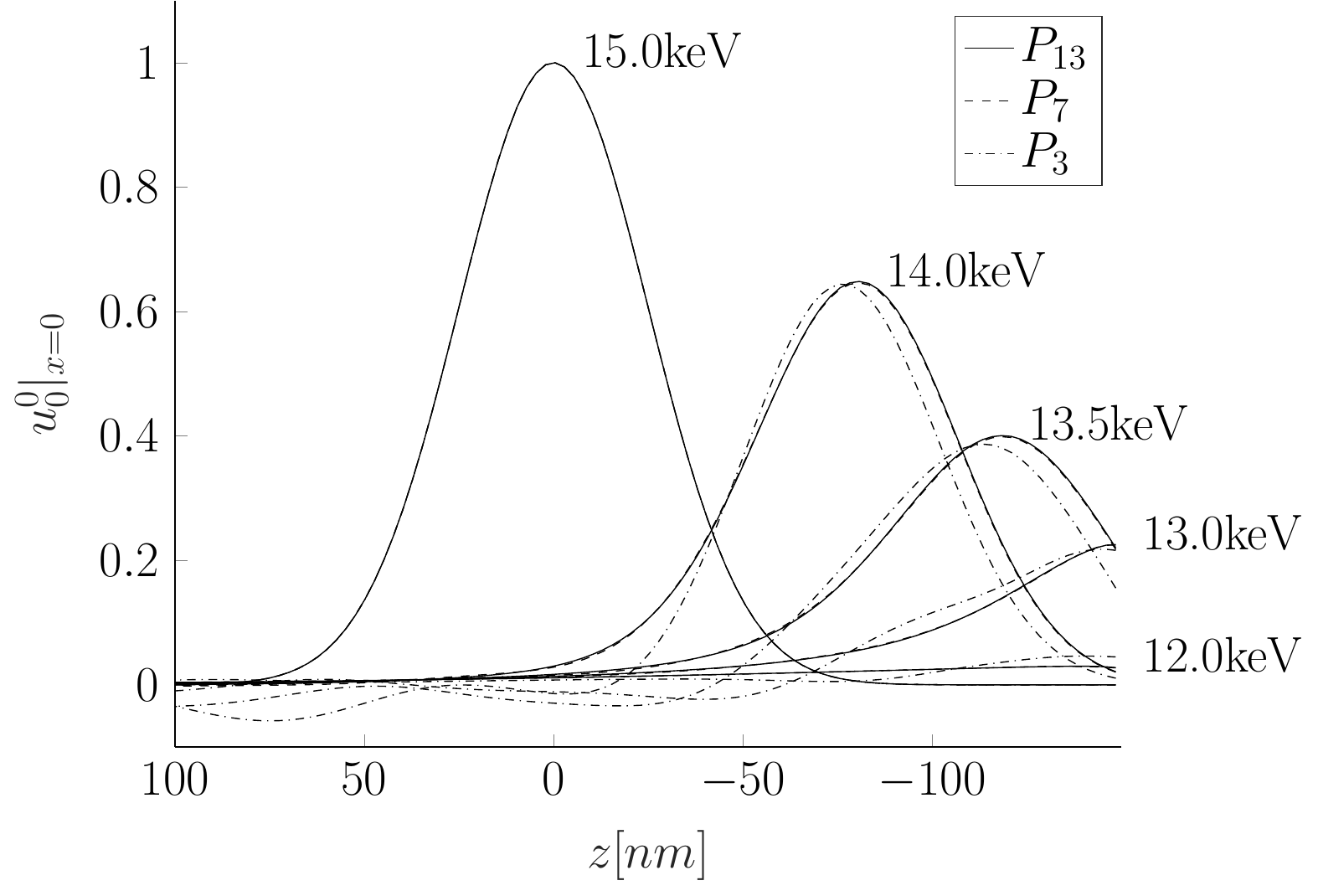}
		\subcaption{\PN{3} and \PN{7} vs. \PN{13}}
		\label{fig:EPMAvacuumOneD2}
	\end{subfigure}
	\caption{Fluence density $u_0^0$ at $x$=0nm of a bulk of electrons inside a material sample of pure copper, initially ($15$keV) moving in negative z-direction, escaping into vacuum over a polished surface at $z=-150$nm computed with \PN{N} and MC.}	
	\label{fig:EPMAvacuum1D}
\end{figure}

First we look at the situation of beam electrons leaving a material sample into vacuum. \Cref{fig:EPMAvacuum2D} shows a \PN{13} result of the electron fluence density $u_0^0$ for a bulk of electrons inside copper, which is initially ($\epsilon=15$keV) centered around $x=z=0$nm and moving in negative $z$ direction
\begin{equation*}
	\psi_{15\text{keV}}(x,z,\Omega) = e^{-\left( \frac{x^2+z^2}{2 \sigma^2 } \right)} \delta(\Omega_z+1)
	\qquad\text{with } \sigma = 25\text{nm},
\end{equation*}
escaping into vacuum over a polished surface at $z=-150$nm. To assess the quality of the \PN{13} result \Cref{fig:EPMAvacuum2D} also shows MC results. For all four snapshots we can observe a very good agreement between \PN{13} and MC, in fact, it is hard to see deviations between the results of the two methods. The good agreement is more visible in \Cref{fig:EPMAvacuumOneD1}, which shows snapshots of $u_0^0$ sliced at $x=0$. \Cref{fig:EPMAvacuumOneD2} shows the same slices obtained with approximation orders $N=3,7,13$. Taking the \PN{13} as reference, we observe that the \PN{7} result is already very accurate whereas the \PN{3} result deviates clearly. Comparing the deviation at the boundary with the deviations in the interior, it seems that the deviation is caused by the low order rather then a poor physical model at the boundary.

\subsubsection{Test Case 4: Gaussian Beam}

\begin{figure}
	\centering
	\begin{subfigure}{0.49\hsize}
		\centering
		\includegraphics[width=0.9\hsize]{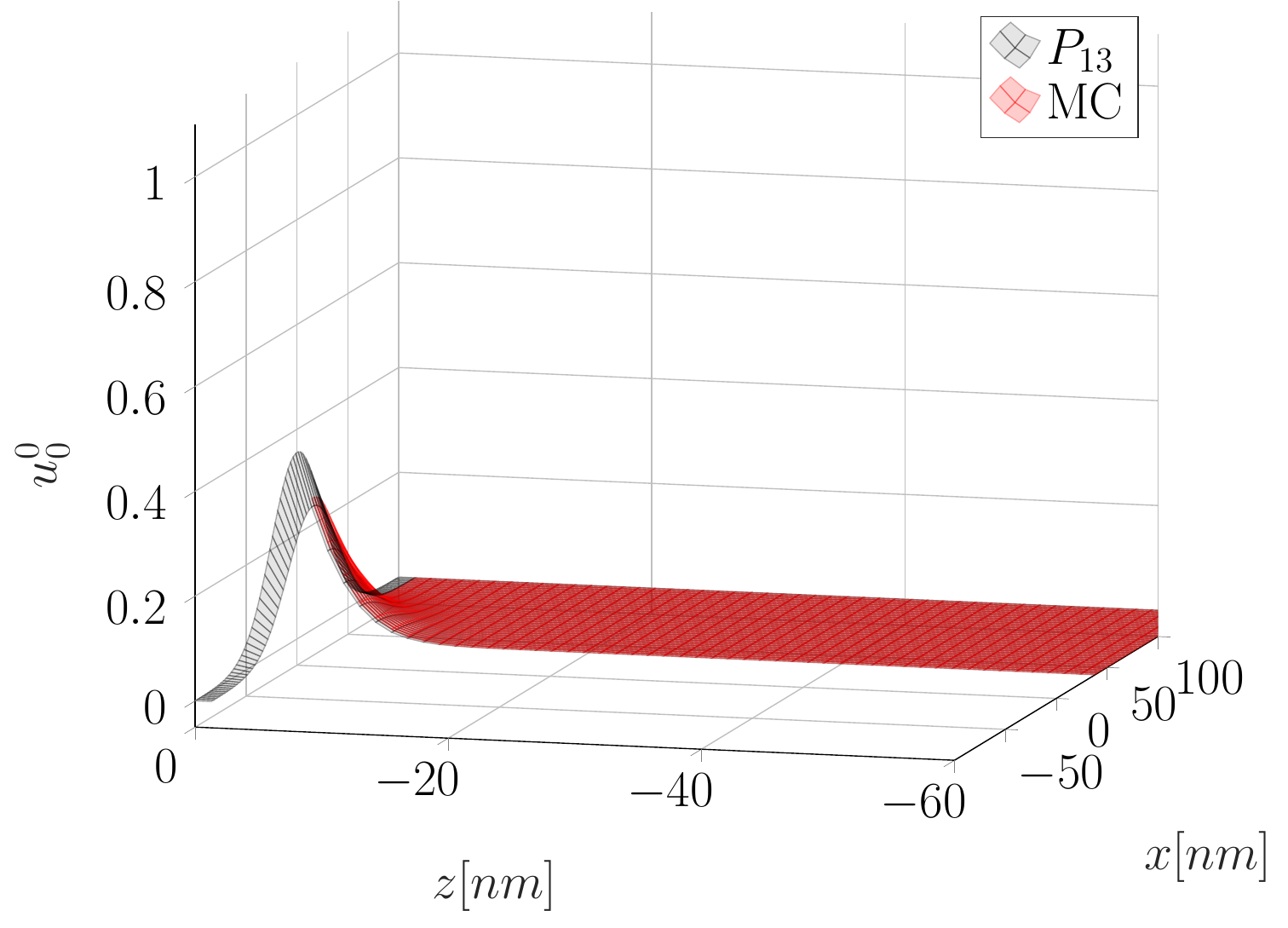}
		\subcaption{$u_0^0$ at $14.1$keV}
	\end{subfigure}
	\begin{subfigure}{0.49\hsize}
		\centering
		\includegraphics[width=0.9\hsize]{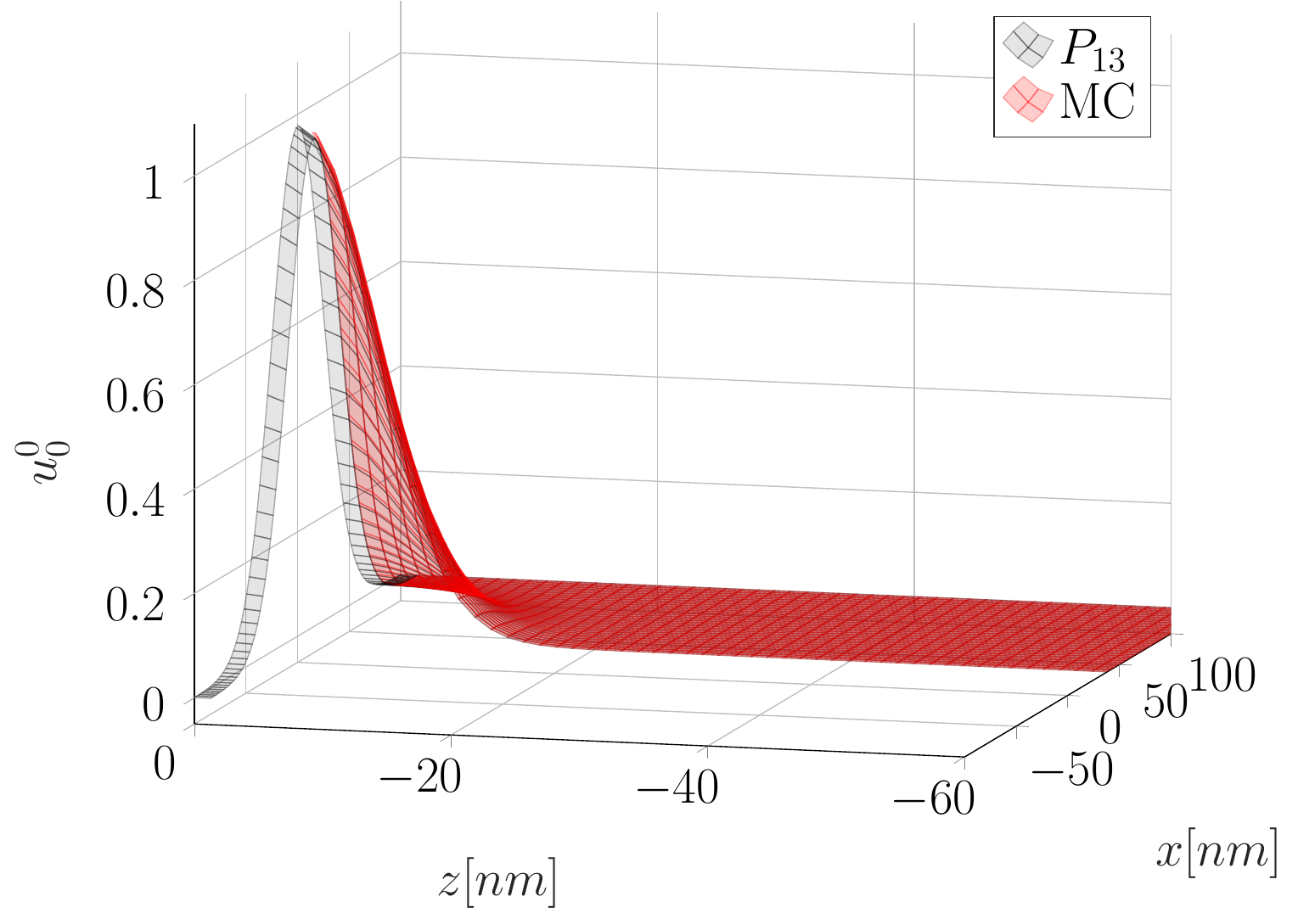}
		\subcaption{$u_0^0$ at $14.0$keV}
	\end{subfigure}
	\begin{subfigure}{0.49\hsize}
		\centering
		\includegraphics[width=0.9\hsize]{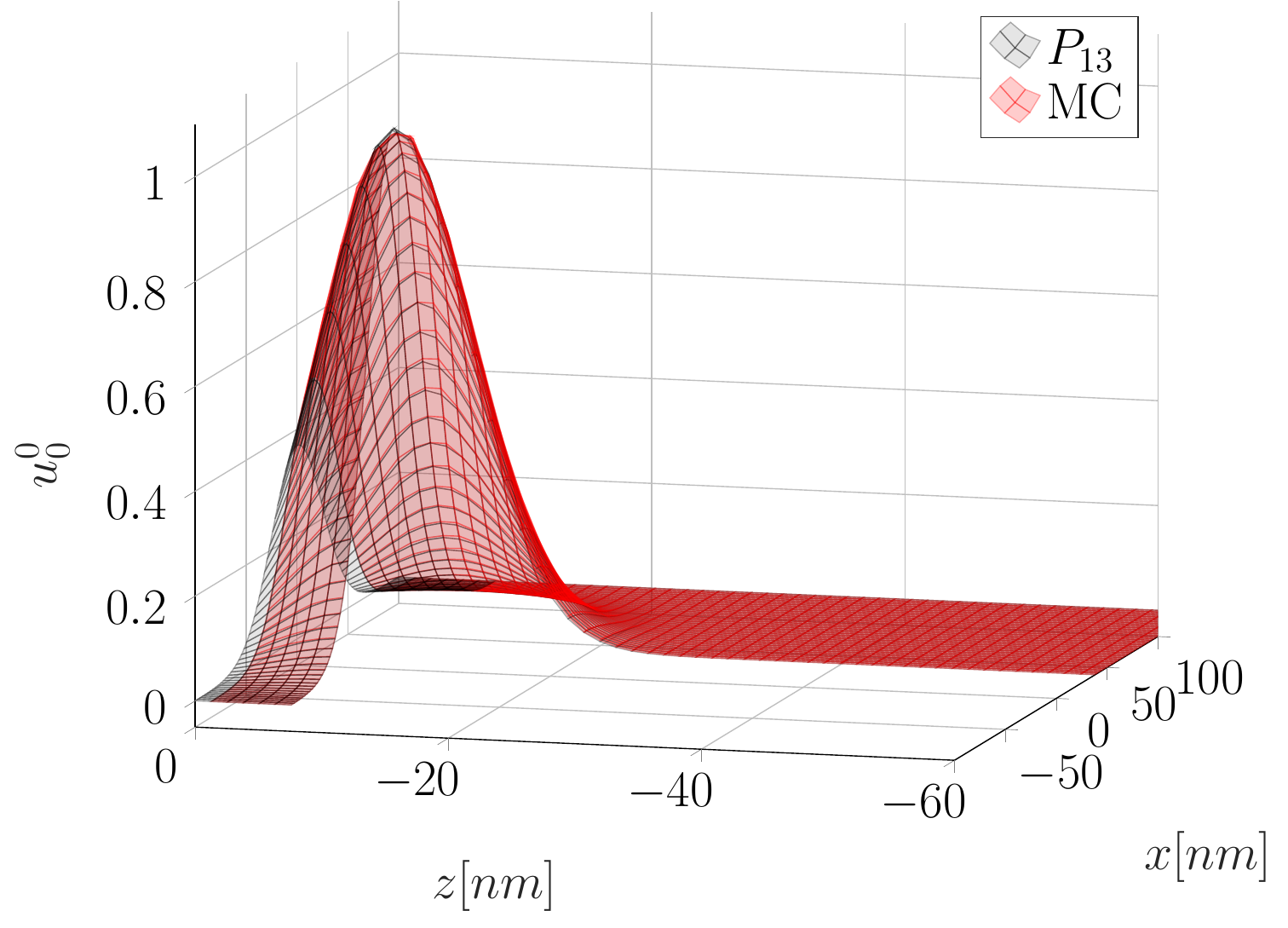}
		\subcaption{$u_0^0$ at $13.9$keV}
	\end{subfigure}
	\begin{subfigure}{0.49\hsize}
		\centering
		\includegraphics[width=0.9\hsize]{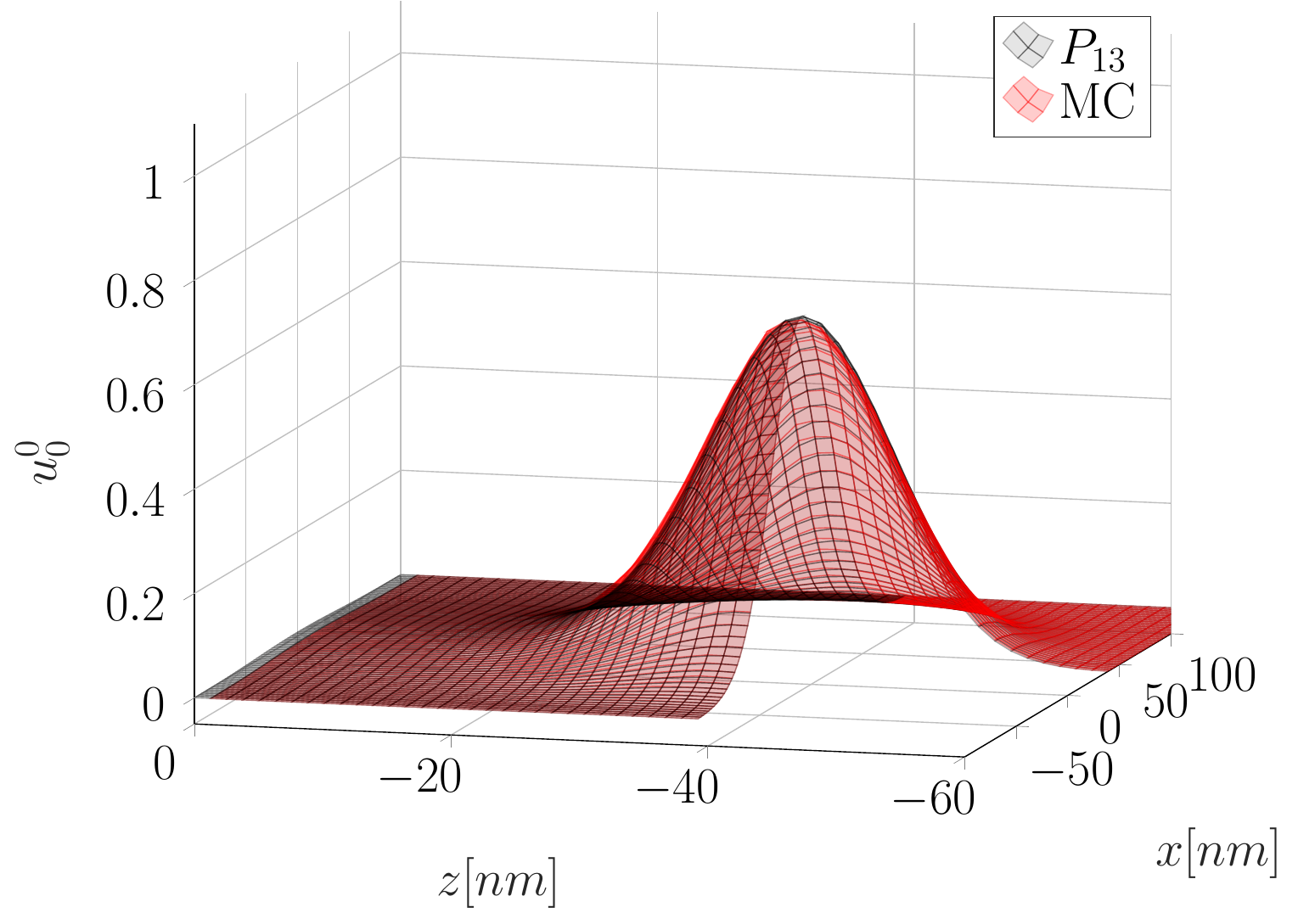}
		\subcaption{$u_0^0$ at $13.5$keV}
	\end{subfigure}
	\caption{Fluence density $u_0^0$ of beam electrons, shot from a 14keV Gaussian electron beam, penetrating a material sample of pure copper over a polished surface at $z$=0nm computed with \PN{13} and MC.}
	\label{fig:gaussianBeamTwoDim}
\end{figure}

\begin{figure}
	\centering
	\begin{subfigure}{0.49\hsize}
		\centering
		\includegraphics[width=0.9\hsize]{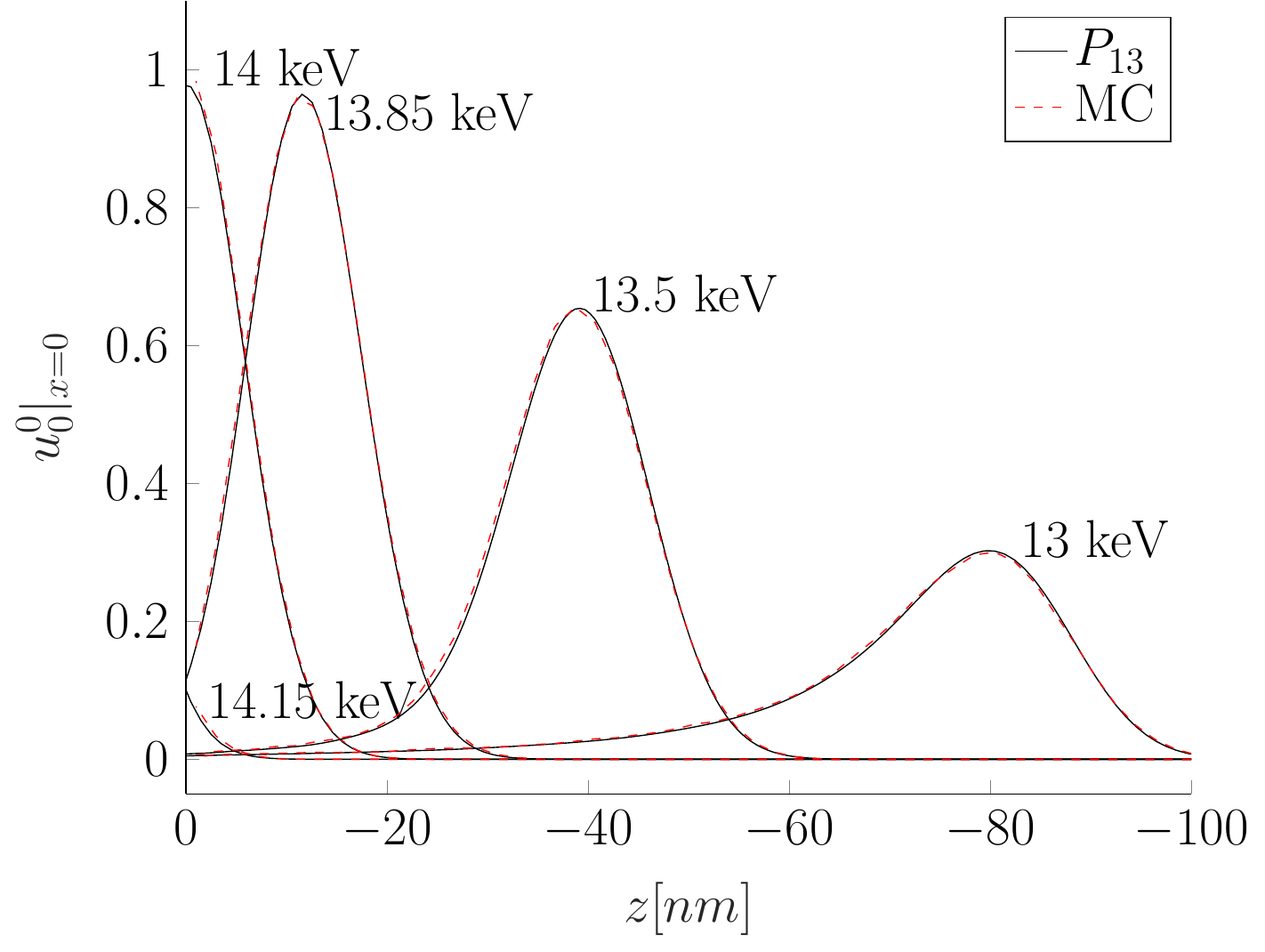}
		\subcaption{\PN{13} vs. MC}
		\label{fig:gaussianBeamOneDimA}
	\end{subfigure}
	\begin{subfigure}{0.49\hsize}
		\centering
		\includegraphics[width=0.9\hsize]{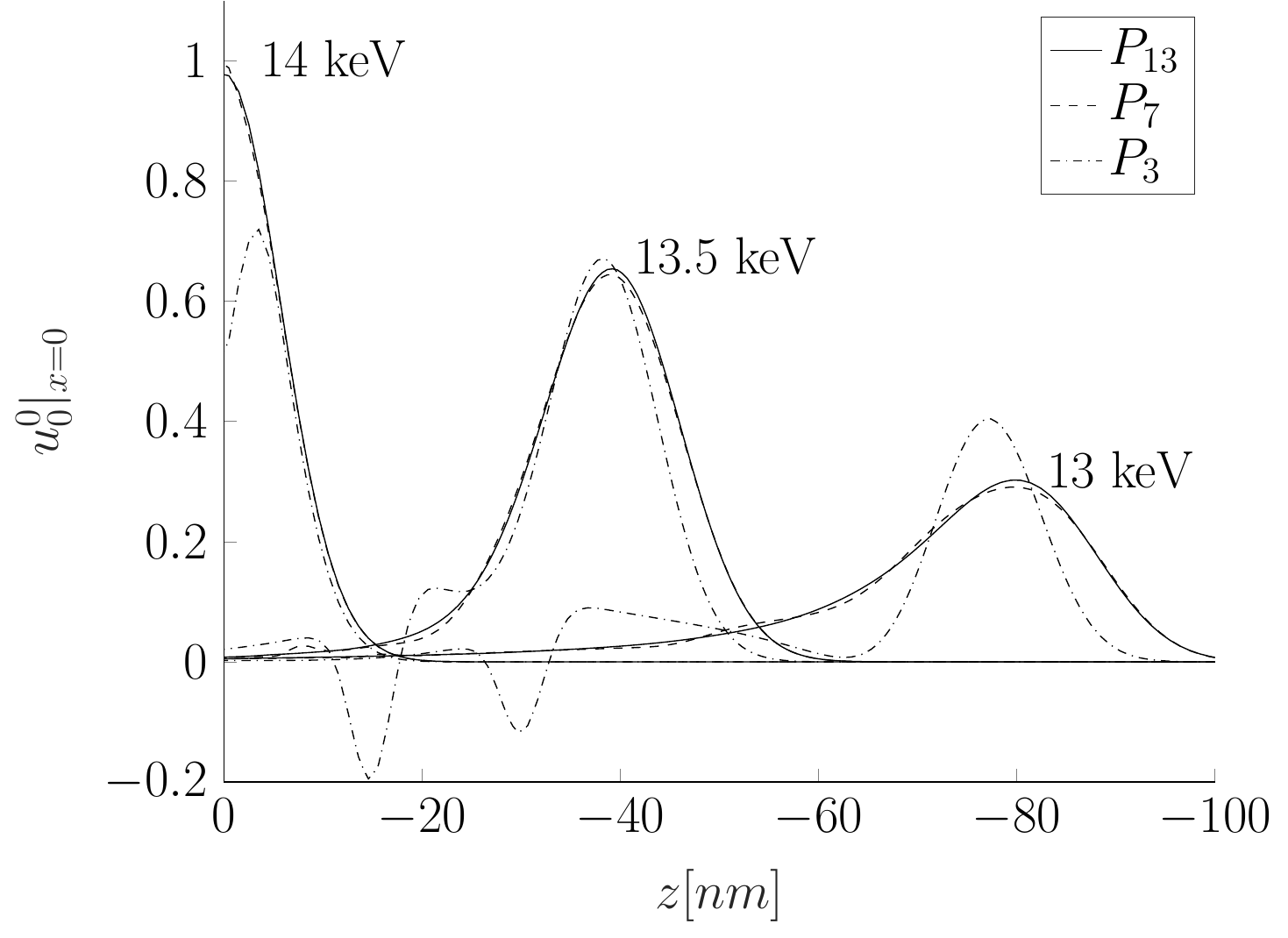}
		\subcaption{\PN{3} and \PN{7} vs. \PN{13}}
		\label{fig:gaussianBeamOneDimB}
	\end{subfigure}
	\caption{Fluence density $u_0^0$ at $x$=0nm of beam electrons, shot from a 14keV Gaussian electron beam, penetrating a material sample of pure copper over a polished surface at $z$=0nm computed with \PN{N} and MC.}
	\label{fig:gaussianBeamOneDim}
\end{figure}

In this test case we look at the other important boundary scenario in the context of EPMA: Beam electrons entering a material probe normal to the polished material surface. As for the previous test case we assume solutions to be constant in the $y$ $(\hat{=}x_2)$ direction and look at a two dimensional scenario in $x$ ($\hat{=}x_1$) and $z$ ($\hat{=}x_3$) direction.\\
We consider a 14keV Gaussian electron beam that is aligned with the z-axis and shoots electrons onto a polished surface located at $z=0$, which we can model via the incoming electron fluence
\begin{equation*}
	\psi_{in}(\epsilon,x,\Omega) =
	e^{-\big( \frac{\epsilon-14\text{keV}}{\sqrt2 \sigma_\epsilon} \big)^2} e^{-\big( \frac{x}{\sqrt2 \sigma_x} \big)^2} e^{-\big( \frac{\Omega_z+1}{\sqrt2 \sigma_\Omega} \big)^2}.
\end{equation*}
In \Cref{fig:gaussianBeamTwoDim} we compare the electron fluence density $u_0^0$ obtained with \PN{13} to a MC result for $\sigma_\epsilon=\frac{14\text{keV}}{100}$, $\sigma_x=25$nm and $\sigma_{\Omega_{z}}=0.1$. As before, \PN{13} reproduces MC very well and no significant deviations between the results can be observed. The good agreement is more apparent in \Cref{fig:gaussianBeamOneDimA} where slices of $u_0^0$  at different energies can be compared. In \Cref{fig:gaussianBeamOneDimB} we compare different approximation orders N and make the same observation as for the test case on vacuum boundary conditions: The \PN{7} result shows only minor deviations from the \PN{13} result, whereas the \PN{3} result deviates significantly. Again, it seems that the strong deviation of \PN{3} from \PN{13} is related to the low approximation order rather then introduced by a poor translation of the Gaussian electron beam into boundary conditions for the \PN{N} equations.

	\section{Acknowledgements}
	The authors acknowledge funding of the German Research Foundation (DFG) under grant TO 414/4-1.
	\newpage
	\bibliography{bibliography} 
	\bibliographystyle{ieeetr}

\end{document}